\documentclass[opre,nonblindrev,12pt]{informs3_modified}
\MANUSCRIPTNO{}       
\JOURNAL{}
\def\theARTICLETOP{}

\OneAndAHalfSpacedXI 

\usepackage{natbib}
 \bibpunct[, ]{(}{)}{,}{a}{}{,}%
 \def\bibfont{\small}%

\usepackage{bbm}
\usepackage{rotating}
\usepackage{booktabs}
\usepackage{amsmath, amssymb}
\usepackage{amsfonts}
\usepackage{mathrsfs}
\usepackage{eurosym}
\usepackage[mathscr]{euscript}
\usepackage{subfigure}

\usepackage{graphicx}

\usepackage[mathscr]{euscript}
\DeclareSymbolFont{rsfs}{U}{rsfs}{m}{n}
\DeclareSymbolFontAlphabet{\mathscrsfs}{rsfs}
\usepackage{float}
\usepackage{soul} 
\usepackage{tikz}
 \usepackage{lmodern}
 \tikzset{
    invisible/.style={opacity=0,text opacity=0},
    visible on/.style={alt={#1{}{invisible}}},
    alt/.code args={<#1>#2#3}{
      \alt<#1>{\pgfkeysalso{#2}}{\pgfkeysalso{#3}} %
      },
   }
\usetikzlibrary{patterns}
\usepackage{pgfplots}
\pgfplotsset{compat=newest}
\usetikzlibrary{arrows,decorations.pathreplacing,backgrounds,automata}
\usepackage{multirow, multicol}
\usepackage{enumerate, geometry}
\usepackage[colorlinks=true, allcolors=blue]{hyperref}

\newcommand{\bA}{ \mathbf{A} }
\newcommand{\ba}{ \mathbf{a} }
\newcommand{\bb}{ \mathbf{b} }
\newcommand{\bbb}{ \mathbf{b} }

\newcommand{\bc}{ \mathbf{c} }
\newcommand{\bD}{ \mathbf{D} }
\newcommand{\bd}{ \mathbf{d} }

\newcommand{\bg}{ \mathbf{g} }
\newcommand{\bG}{ \mathbf{G} }

\newcommand{\bh}{ \mathbf{h} }

\newcommand{\bP}{ \mathscrsfs{D} }
\newcommand{\mf}{ \mathscr{F} }

\newcommand{\bw}{ \mathbf{w} }

\newcommand{\bx}{ \mathbf{x} }

\newcommand{\by}{ \mathbf{y} }

\newcommand{\bzero}{ \mathbf{0} }

\newcommand{\cU}{\mathcal{U}}
\newcommand{\cA}{\mathcal{A}}
\newcommand{\cB}{\mathcal{B}}
\newcommand{\cC}{\mathcal{C}}
\newcommand{\cI}{\mathcal{I}_1}
\newcommand{\cG}{\mathcal{G}}
\newcommand{\cS}{\mathcal{S}}
\newcommand{\cJ}{\mathcal{J}}
\newcommand{\cX}{\mathcal{X}}

\newcommand{\cK}{\mathcal{K}}
\newcommand{\cH}{\mathcal{H}}
\newcommand{\hatA}{\hat{\bA}}
\newcommand{\hatb}{\hat{\bb}}
\newcommand{\barA}{\Delta} 
\newcommand{\barhatA}{\hat{\Delta}}
\newcommand{\bone}{\mathbf{1}}

\newcommand{\FO}{\mathbf{FO}}
\newcommand{\IO}{\mathbf{IO}} 
\newcommand{\MIO}{\mathbf{MIO}} 
\newcommand{\TMIO}{\mathbf{MIO}}

\newcommand{\eMIO}{\mathbf{e\text{-}MIO}} 
\newcommand{\zeroObj}{Indifference Measure}
\newcommand{\priorB}{Adherence Measure}
\newcommand{\minDist}{Adjacency Measure}
\newcommand{\fairness}{Fairness Measure}
\newcommand{\minmin}{Compactness Measure}

\newcommand\hatX{\mathbf{x}^0}

\newcommand{\ie}{{\it i.e.}}
\newcommand{\eg}{{\it e.g.}}

\definecolor{myGreen}{RGB}{0,0,0}
\definecolor{myPink}{RGB}{255,153,153}
\definecolor{myCyan}{RGB}{0,206,209}

\TheoremsNumberedThrough     

\EquationsNumberedThrough    

\def\add#1{{\color{black} #1}}

\newcolumntype{L}[1]{>{\raggedright\let\newline\\arraybackslash\hspace{0pt}}m{#1}}
\newcolumntype{C}[1]{>{\centering\let\newline\\arraybackslash\hspace{0pt}}m{#1}}
\newcolumntype{R}[1]{>{\raggedleft\let\newline\\arraybackslash\hspace{0pt}}m{#1}}

\begin{document}

\TITLE{Inferring Linear Feasible Regions \\using Inverse Optimization\footnote{Both authors contributed equally to this manuscript.}
}
\RUNTITLE{Inferring Linear Feasible Regions using Inverse Optimization}

\RUNAUTHOR{Ghobadi and Mahmoudzadeh}

\ARTICLEAUTHORS{
\AUTHOR{Kimia Ghobadi}
\AFF{Malone Center for Engineering in Healthcare, Center for Systems Science and Engineering, Department of Civil\\ and Systems Engineering, Johns Hopkins University, Baltimore, MD, USA. \EMAIL{kimia@jhu.edu}}
\AUTHOR{Houra Mahmoudzadeh}
\AFF{Department of Management Sciences, University of Waterloo, ON, Canada. \EMAIL{houra.mahmoudzadeh@uwaterloo.ca}}
}

\ABSTRACT{%
Consider a problem where a set of feasible observations are provided by an expert and a cost function is defined that characterizes which of the observations dominate the others and are hence, preferred. Our goal is to find a set of linear constraints that would render all the given observations feasible while making the preferred ones optimal for the cost (objective) function. \add{By doing so, we infer the implicit feasible region of the linear programming problem}. Providing such feasible regions ({\it i}\,) builds a baseline for categorizing future observations as feasible or infeasible, and ({\it ii}\,) allows for using sensitivity analysis to discern changes in optimal solutions if the objective function changes in the future. 

In this paper, we propose an inverse optimization framework to recover the constraints of a forward optimization problem using multiple past observations as input. We focus on linear models in which the objective function is known but the constraint matrix is partially or fully unknown. \add{We propose a general inverse optimization methodology that recovers the complete constraint matrix and then introduce a tractable equivalent reformulation.} Furthermore, we provide and discuss several generalized loss functions to inform the desirable properties of the feasible region based on user preference and historical data. 
Our numerical \add{examples} verify the validity of our approach, emphasize the differences among the proposed measures, and provide intuition for large-scale implementations. \add{We further demonstrate our approach using a diet recommendation problem to show how the proposed models can help impute personalized constraints for each dieter.}} 


\maketitle

\section{Introduction} \label{sec:introduction}
Conventional (forward) optimization problems find an optimal solution for a given set of parameters. Inverse optimization, on the other hand, infers the parameters of a forward optimization problem given a set of observed solutions (typically a single point). In the literature, inverse  optimization \citep{zhang1996calculating} is often employed to derive the parameters of the cost vector of an optimization problem while the constraint parameters are assumed to be fully known. In this paper, we focus on the opposite case. We impute the constraint parameters (as opposed to the objective function) of a linear forward problem given a cost vector and a set of observations. \add{Hence, we infer the feasible region of the forward problem which can be used to identify future feasible or infeasible observations, and to understand the behaviour of the model under different cost vectors.}

When imputing the cost vector, it is usually assumed that the observed solution is a candidate optimal solution~\citep{Ahuja01, Iyengar05, ghate2020inverse}. More recently, several studies also consider the case where the observed solution is not necessarily a candidate for optimality. They propose inverse models to minimize error metrics that capture the optimality gap of the observed solution~\citep{Keshavarz11, Chan14, Chan15, Bertsimas15, Aswani16, naghavi2019inverse}.  More recently, multiple observations are considered, where the cost vector is imputed based on a given set of feasible  observations~\citep{Keshavarz11,Troutt06,Troutt08,Chow12, Bertsimas15,  esfahani2018IncompleteInfo, babier2019ensemble}. 
\citet{tavasliouglu2018structure} find a set of inverse-feasible cost vectors, instead of a single cost vector, that makes feasible observations optimal. \add{A standard assumption in the literature of inverse optimization is that the observed data is noise-free \citep{zhang1999further, Ahuja01}. There are only a few studies that consider noise or uncertainty in the input data \citep{Aswani16, dong2018inferring, ghobadi2018robust} when inferring the cost vector.}

Extending from only imputing the cost vector, some studies consider the case where both the objective function and the right-hand side (RHS) of the constraints are imputed simultaneously for specific types of problems~\citep{dempe2006inverse, Chow12, cerny2016inverse}. Note that when the feasible region is being imputed, a given observation can become optimal since the constraints can be adjusted so as to position the given observation on the boundary of the feasible region. A few studies focus on imputing only the RHS constraint parameters of the forward problem. Given a single observation, \citet{cerny2016inverse} find the RHS of the constraints from a pre-specified set of possible parameters. In other studies, the RHS is imputed in such a way that the observed solution becomes optimal~\citep{birge2017inverse} or near-optimal according to a pre-specified distance metric~\citep{dempe2006inverse, guler2010capacity, saez2018short}. 

\add{There are limited previous studies that impute the left-hand side (LHS) parameters of the constraints set. \citet{chan2018inverse} input a single observation and propose an inverse optimization method to find the LHS (assuming the RHS is known) as well as the cost vector such that the given observation becomes optimal. Assuming an unknown objective in addition to an unknown feasible region would result in finding a feasible region that makes a given point optimal for \emph{some} objective (\ie, any objective that fits the problem mathematically), and hence, would make the problem more relaxed and less practical. An unknown objective would also mean that we are not able to assess the quality of the given solutions.}

\add{Although the forward problem we consider is similar to that of~\cite{chan2019inverse}, our proposed inverse models differ from theirs in several key aspects. Our models infer the full constraint parameters (both LHS and RHS) and consider multiple observations instead of a single one. We assume that the objective function is known and hence, we can identify the preferred solution(s) among all given observations. This assumption is relevant in practical settings and will be discussed in Section~\ref{sec:motivation}. Their models also require a prior belief and additional user-defined conditions on the constraint parameters  to avoid trivial (all zero) solutions; our models do not require any additional user input and generate non-trivial solutions by design. We also introduce generalized loss functions that do not necessarily rely on a prior belief on the constraint parameters.}

Solving inverse optimization problems efficiently has been the focus of a few papers in the literature. While inverse optimization problems for imputing the cost vector often retain the complexity of their corresponding forward problems (\eg, linear programming), imputing the constraint parameters often constitutes a non-convex problem due to the presence of multiple bilinear terms. Hence, the resulting inverse problems are more complex to solve. To address these concerns, a few studies in the literature focus on specific problem instances and find certain criteria or assumptions to reduce the complexity~\citep{birge2017inverse, brucker2009inverse}. Others propose a solution methodology that solves a sequence of convex optimization problems under a specific distance metric~\citep{chan2018inverse}. In this paper, instead of attempting to solve a nonlinear non-convex problem, we use the problem's theoretical properties and propose an equivalent reformulation that can be linearly constrained, and hence, easier to solve. We also further simplify the problem by providing closed-form solutions or suggesting decomposition approaches for specific cases. 

To the best of our knowledge, this paper is the first in the literature to propose an inverse optimization framework for inferring the full constraint matrix (both LHS and RHS) of a linear programming model based on multiple observations.  \add{Contrary to the recent literature, the objective function is known in our proposed framework, and the constraint parameters are partially or fully unknown. Our models do not require any additional input data with the observations, but such data can be incorporated in the model if available. The solutions are observed without noise, but our framework allows for potential inclusion of noisy data.} 
The contributions of this paper are \add{summarized} as follows: 
\begin{itemize}
    \item We propose a single-point inverse optimization model that inputs one observation and infers a set of fully or partially unknown constraint parameters of the forward problem so as to make the given observation optimal for a specific cost vector.
    \item We extend this model to a multi-point inverse optimization methodology that inputs any number of observations and finds the constraint parameters so as to make all the observations feasible and the preferred observation(s) optimal for a specific cost vector. 
    \item We develop an equivalent tractable reformulation of the multi-point inverse optimization model that eliminates the bilinearity of the original model. 
    \item \add{We introduce a generalized loss function that can take any form or input any type of available data and induces the desirable properties of the feasible region of the forward problem. This proposed loss function does not necessarily rely on a prior belief, expert opinion, or other specific user inputs on the constraint parameters.}
    \item We test and validate our proposed methodology on numerical \add{examples} and demonstrate the characteristics of each of the loss functions introduced. 
    \item \add{We demonstrate the application of our approach on a diet recommendation problem and show that the proposed model can use past food consumption observations to impute each user's implicit constraints and generate personalized diets that are palatable.} 
\end{itemize}

The rest of this paper is organized as follows. In  Section~\ref{sec:motivation}, we motivate  the proposed methodology by presenting examples of application areas. Section~\ref{sec:Methodology} introduces our methodology for inverse optimization of constraint parameters, its theoretical properties, and an equivalent reformulation. In Section~\ref{sec:measures}, we introduce examples of the generalized loss function that can be used as the objective function of the inverse optimization problem and discuss the theoretical properties of each. We illustrate the results of our methodology using two numerical examples and a diet recommendation application in Section~\ref{sec:numericalexample}. Section~\ref{sec:discussions} discusses a few extensions to the proposed models, and finally, conclusions and future research directions are provided in Section~\ref{sec:conclusions}.

\section{Motivation} \label{sec:motivation}
Inverse optimization has been applied to several different application areas, from healthcare~\citep{Erkin10, Ayer14} and nutrition~\citep{ghobadi2018robust} to finance~\citep{Bertsimas12} and electricity markets~\citep{birge2017inverse}. In this section, we provide two example applications where imputing the feasible region based on a set of collected observations is of practical importance. These applications serve to motivate the development of our proposed methodology. 

\vspace{1em}
\noindent{\bf A. Cancer Treatment Planning:} 
Consider the radiation therapy treatment planning problem for cancer patients. The input of the problem is a medical image (\eg, CT, MRI) which includes contours that delineate the cancerous region (\ie, tumor) and the surrounding healthy organs. The goal of a treatment planner is to find the direction, shape, and intensity of radiation beams such that a set of clinical metrics on the tumor and the surrounding healthy organs are satisfied.  \add{While there exists literature on using inverse optimization for inferring the objective function in cancer treatment planning \citep{Chan14, Goli15, babier2018inverse}, to the best of our knowledge, no study infers the constraint parameters for cancer treatment planning.}

In current practice, there are clinical guidelines on the upper/lower limits for different clinical metrics. Planners often try to find an \emph{acceptable} treatment plan based on these guidelines to forward to an oncologist who will inspect it and either approve or return it to the planner. If the plan is not approved, the planner receives a set of instructions on which metrics to adjust. It often happens that the final approved plan may not meet all the clinical limits simultaneously as there are trade-offs between different metrics.  

Suppose we have a set of approved treatment plans from previous patients. Even though there are clinical guidelines on acceptability thresholds for different metrics, in reality, there may exist approved treatment plans that do not meet these limits. There may also exist plans that do meet the guidelines but are not approved by the oncologists since she/he believed a better plan is achievable. Hence, the true feasible region of the forward problem in treatment planning is unknown. 
Considering the historically-approved plans as ``feasible points'', we can employ our inverse optimization approach to find the constraint parameters, based on which we can understand the implicit logic of the oncologists in approving a treatment plan. In doing so, we would help both the oncologists and the planners by ({\it i}\,) generating more realistic lower/upper bounds on the clinical metrics based on past observations, ({\it ii}\,) improving the iterative planning process by producing higher quality initial plans given the clear guidelines and hence, reducing the number of preliminary plans passed back and forth between the planner and oncologists, and ({\it iii}\,) improving the quality of plans by preventing low-quality solutions that otherwise satisfy the acceptability metrics, especially for inverse treatment planning methods that heavily rely on these metrics.

\vspace{1em}
\noindent{\bf B.  Diet Problem:} 
Consider a diet recommendation system that suggests a variety of food items based on a user's dietary needs and/or personal preferences. Each meal can be characterized by a set of features and/or metrics such as meat content or daily value of each nutrient.The diet problem often consists of minimizing some objective such that a set of requirements on the food intake is met.  \add{A limited number of studies in the literature have used inverse optimization for inferring the objective function weights in a diet problem \citep{shahmoradi2019quantile, ghobadi2018robust}. However, to the best of our knowledge, no study infers the constraint parameters for the diet problem.}

Assume the objective function of the underlying (forward) optimization problem in the diet problem is known. Examples of such objective functions would be minimizing total calories in a weight loss program, minimizing sodium intake in a hypertension diet, or minimizing monetary cost. In addition to dietary requirements, each person has a set of implicit constraints that would result in them finding a certain suggestion ``palatable'' or not. Different users would have different such constraints, and it is not explicitly possible to list what these constraints are. In such cases, our inverse optimization model can use historical data to ensure the next suggested meal in the diet recommendation system is palatable. 

For example, consider a user who is mostly vegetarian and is implicitly limiting the number of meat servings during the week, or another user who prefers to limit the amount of dairy intake when consuming certain vegetables. If we observe the user's diet choices (feasible observations) for a certain time horizon, the inverse optimization model would find the set of constraints (feasible region) that captures this behaviour by making diets that are too far off from the past observations infeasible while ensuring that the required amount of nutrients are met, the diet is palatable (feasible), and the given objective (\eg, cost, calories) is minimized.  \add{In Section~\ref{sec:numericalexample}, we further discuss this application of our proposed methodology on a diet recommendation problem. }

\section{Methodology} \label{sec:Methodology}
In this section, we first set up the forward optimization problem where, contrary to conventional inverse optimization, the cost vector is known and the unknown parameters are, instead, the constraint parameters. Let $\bc \in\mathbb{R}^n,\bA\in\mathbb{R}^{m_1\times n},\bbb\in\mathbb{R}^{m_1}, \bG\in\mathbb{R}^{m_2\times n}$ and $\bh\in\mathbb{R}^{m_2}$. We define our linear forward optimization ($\FO$) problem as 
\begin{subequations}
\begin{align}
\FO: \quad \underset{\bx}{\text{minimize}} & \quad \bc'\bx\\
\text{subject to} & \quad \bA\bx \ge \bbb,  
\\ 
& \quad \bG \bx \ge \bh,  \\
&{\color{myGreen} \quad \bx \in \mathbb{R}^n. }
\end{align}
\end{subequations}

Consider the case in which some or all of the constraint parameters of the $\FO$ formulation are unknown, but there exist one or more observations (solutions) that are deemed feasible or optimal for $\FO$ based on expert opinion. 
For such settings, we propose inverse optimization models that infer these unknown parameters of $\FO$ and recover the full feasible region. We assume that $\bA$ and $\bb$ are the \emph{unknown constraint} parameters that the inverse optimization aims to infer and $\bG$ and $\bh$ are the \emph{known constraint} parameters.

For every constraint of $\FO$, three cases can be considered: ({\it i}\,) all of the parameters of the constraint are known, in which case we denote it as part of the known constraints, $\bG \bx \geq \bh$; ({\it ii}\,) all of its parameters are unknown and we denote the constraint as part of the unknown constraints, $\bA \bx \geq \bb$; or ({\it iii}\,) some of its parameters are known (for instance, $\bb$ is known) or some properties about the constraint(s) are known, in which case we denote the constraint(s) as part of unknown constraint parameters $\bA$ and $\bb$ and add the additional restrictions to the inverse model. We will discuss the latter case in more detail in Section~\ref{sec:discussions}, and without loss of generality, we assume no such partial information is available for the rest of this section.

Throughout this paper, we index the unknown and known constraints by the sets $\cI =\{ 1,\dots, m_1\}$ and $\mathcal{I}_2 =\{ 1,\dots, m_2\}$, respectively. Note that if there are no known constraints, $\mathcal{I}_2 = \emptyset$. The $i^\text{th}$ row of the constraint matrices $\bA$ and $\bG$ is referred to as $\ba_i$ and $\bg_i$, respectively. Similarly, the $i^\text{th}$ elements of the $\bb$ and $\bh$ vectors are denoted by $b_i$ and $h_i$, respectively.  The set $\cJ= \{1,\dots,n\}$ denotes the columns in the constraint matrices (\ie, the indices of the $\bx$ variable). We use bold numbers $\mathbf{1}$ and $\bzero$ to denote the all-ones and the all-zeros vectors, respectively. 

In this section, we propose three models to infer the unknown parameters of $\FO$. First, we present a single-point inverse optimization model when only one observation is available. Next, we provide a multi-point inverse optimization model to infer the unknown constraint parameters when several observations are available. Finally, we propose a tractable reformulation for the proposed inverse optimization model.

\subsection{Single-point Inverse Optimization}
In this section, we propose an inverse optimization model for the case where only one observed solution is available. Given a single observation $\hatX$, a cost vector $\bc$, and known constraint parameters $\bG$ and $\bh$ (if any), we would like to formulate an inverse optimization model that finds the unknown constraint parameters $\bA$ and $\bb$ such that the observation $\hatX$ is optimal for the forward problem $\FO$. Without loss of generality, we assume that the observation $\hatX$ is feasible for the known constraints $\bG \bx \geq \bh$, since the forward problem will be otherwise ill-defined,  and the inverse problem will have no solution.

Let $\by$ and $\bw$ be dual vectors for constraints (1b) and (1c) of $\FO$, respectively. The single-point inverse optimization model ($\IO$) can be written as follows: 
\begin{subequations} \label{eq:IO}
\begin{align}
\IO: \underset{\bA, \bb, \by, \bw}{\text{minimize}} & \quad \mf(\bA, \bb; \bP) ,  \\
\text{subject to} 
& \quad \bA \hatX \geq \bbb, \label{eq:IOprimalfeas}\\ 
& \quad \bc' \hatX = \bbb'\by + \bh' \bw,  \label{eq:IOstrongduality}\\ 
& \quad \bA'\by + \bG' \bw = \bc,   \label{eq:IOdualfeas1}\\ 
&\quad ||\ba_i||= 1, \quad \forall i \in \cI  \label{eq:IOnorm} \\
& \quad \by \in \mathbb{R}^{m_1}, \quad \bw \in \mathbb{R}^{m_2}, 
\label{eq:IOdualfeas2} \\
& \quad \bA \in \mathbb{R}^{m_1 \times n}, \quad \bb \in \mathbb{R}^{m_1}. \label{eq:IOprimalfeas2}
\end{align}
\end{subequations}
The objective $\mf(\bA, \bb; \bP)$ is a loss function that drives the desired properties of the feasible region based on some given input parameter $\bP$. \add{For instance, the user may input a prior belief on the shape of the feasible region as parameter $\bP$ and set the objective function $\mf$ to minimize the deviation from such prior belief.} More details and other explicit examples of the loss function $\mf(\bA, \bb; \bP)$ are discussed in Section~\ref{sec:measures}. Constraint~\eqref{eq:IOprimalfeas} enforces primal feasibility of $\bx^0$.  Constraints~\eqref{eq:IOdualfeas1} and~\eqref{eq:IOdualfeas2} are the dual feasibility constraints. Constraint~\eqref{eq:IOstrongduality} is the strong duality constraint that ensures $\bx^0$ is the optimal solution of $\FO$. Finally, without loss of generality, constraint~\eqref{eq:IOnorm} normalizes the LHS of each unknown constraint based on some norm $||\cdot||$. The introduction of this norm avoids finding multiple scalars of the same constraint parameters or finding trivial (all-zero) solutions, without requiring the user to define application-specific side constraints. Nevertheless, if any such side constraints or partial information on $\bA$ or $\bb$ exists, they can be incorporated in the model. This extension will be discussed later in Section~\ref{sec:discussions}.
We make the following assumption to ensure that the forward problem is not a simple feasibility problem. 
\begin{assumption} \label{assum:c}
$\bc \neq  \bzero$.
\end{assumption}
\noindent We note that without Assumption~\ref{assum:c}, the $\IO$ problem will be simplified since it will have many trivial solutions such as $\bA = \bG, \bb= \bh$, and $\by = -\bw$, or alternatively, $\bw=\by=\bzero$ with any $\bA$ and $\bb$ that satisfy the primal feasibility constraint~\eqref{eq:IOprimalfeas}. For the rest of this paper, we assume Assumption~\ref{assum:c} holds. We next show that the $\IO$ formulation is valid and has non-trivial feasible solutions. 
\begin{proposition} \label{prop:IOfeas}
\add{The feasible region of $\IO$ is non-empty.}  
\end{proposition}
\proof{Proof.}
Let $\bw = \bzero$, $\by = || \bc|| \bone$,  \,
$\bb = (\bc'\bx^0)/{||\bc||}$, and $\ba_i= \bc/{||\bc||}$,\, $\forall i \in \cI$, given $\bc\neq \bzero$. 
Then $(\bA, \bb, \by, \bw)$ is a feasible solution to $\IO$. \Halmos
\endproof

\add{In general, any feasible region that renders the point $\bx^0$ optimal for $\FO$ would be a feasible solution to the single-point $\IO$ problem. Therefore, the solutions to single-point $\IO$ can be futile, and the applicability of this model can be limited. For instance, when $m\geq n$, $\IO$ may force all constraints to pass through $\bx^0$ and possibly make $\bx^0$ be the only feasible point for $\FO$. As another example, often more than one feasible observation is available for the forward problem in practice. In such cases, the $\IO$ formulation does not apply because strong duality may not necessarily hold for multiple observations and lead to the infeasibility of some of the observations. As a result, the theoretical properties of the solutions to the inverse optimization problem would also change. Hence, in the next section, we extend our $\IO$ formulation for the case of multiple observations and discuss its properties.}

\subsection{Multi-point Inverse Optimization}

Consider a finite number of observations $\bx^k, k\in \cK =\{1,...,K\}$ to the forward problem. One of the first goals of a multi-point inverse optimization is to find the constraint parameters in such a way that all observations $\bx^k, k\in \cK$ become feasible. We define this property as follows. 

\begin{definition}\label{def:valid}
A polyhedron $\cX = \{\bx\in \mathbb{R}^n \,|\, \bD \bx \geq \bd \}$  is a \underline{valid} feasible set 
for the forward problem if $\bx^k \in \cX$, \, $\forall k\in \cK$. \end{definition}

\begin{remark}\label{rem:validSet}
If $\cX\subseteq \mathbb{R}^n$ is a valid feasible set, then any set $S \subseteq \mathbb{R}^n$ such that $\cX \subseteq S$ is also a valid feasible set. 
\end{remark}

Remark~\ref{rem:validSet} states that if a set $\cX$ is a valid feasible set, then any set that contains $\cX$ is also valid since all observations remain feasible for any superset of $\cX$. Any set that is not valid, \ie, does not contain some observation $\bx^k, k \in \cK$ cannot be a feasible set to the forward optimization (by definition). Hence, all feasible regions that are imputed from the solutions of a multi-point inverse optimization must be valid feasible sets. In particular,  to ensure that the feasible region of the forward problem is well-defined, we assume that the set defined by the known constraints is also a valid feasible set. 
\begin{assumption} \label{assum:G}
The set $\cG = \{\bx \in \mathbb{R}^n |~\bG \bx \geq \bh\}$ is a valid feasible set.
\end{assumption}

\noindent The feasibility of the observations for the known constraints is similar to the assumption in the single-point $\IO$, except that all observations (as opposed to one observation) are assumed to be feasible for the known constraints $\bG \bx \geq \bh$. Otherwise, the inverse optimization will not have a solution.  Although we have $K$ observations in the multi-point inverse optimization, we can identify the observation(s) that result in the best objective function value for the forward problem, because the $\bc$ vector is known. We define the observation with the best value as the preferred observation, denoted by $\bx^0$, for which strong duality must hold. 

\begin{definition}\label{def:x0} The \underline{preferred} solution in a set of observations $\{\bx^k\}_{k\in\cK}$ is defined as
\begin{equation*}
\bx^0 \in \argmin_{\bx^k, k\in\cK}\{\bc'\bx^k\}.    
\end{equation*}
\end{definition}
\noindent 
If multiple observations satisfy Definition~\ref{def:x0}, without loss of generality, we arbitrarily select one of them as $\bx^0$. The multi-point inverse optimization problem aims to find a set of constraints for the forward problem such that all observations are feasible, and the preferred solution $\bx^0$ becomes optimal. 

The following multi-point inverse optimization ($\MIO$) formulation finds a feasible region that minimizes some loss function of the inverse optimal solution (desired properties of the feasible region) from a set of input parameters $\bP$.
\begin{subequations}\label{eq:TMIO}
\begin{align}
\TMIO: \underset{\bA, \bb, \by, \bw}{\text{minimize}} \quad & \mf(\bA, \bb; \bP), \\
\text{subject to} \quad 
& \bA\bx^k \ge \bbb, \quad \forall k\in \cK \label{eq:TMIO-PrimalFeas}\\ & \bc'\bx^0 = \bbb'\by + \bh'\bw, \label{eq:TMIO-StrongDual}\\ 
& \bA'\by  + \bG' \bw = \bc, \label{eq:TMIO-DualFeas}\\ 
& ||\ba_i||= 1, \quad \forall i \in \cI  \label{eq:MIOnorm} \\
& {\color{myGreen}   \by \in \mathbb{R}^{m_1}, \quad \bw \in \mathbb{R}^{m_2},  
\label{eq:MIO-Dualfeas2}} \\
& {\color{myGreen} \bA \in \mathbb{R}^{m_1 \times n}, \quad \bb \in \mathbb{R}^{m_1}. }\label{eq:MIO-signs}
\end{align}
\end{subequations} 
The constraints in $\MIO$ include strong duality~\eqref{eq:TMIO-StrongDual}, dual feasibility (\eqref{eq:TMIO-DualFeas} and~\eqref{eq:MIO-Dualfeas2}), and normalization~\eqref{eq:MIOnorm}. \add{We note that even though $\bc$ is known, the strong duality constraint does not automatically hold. To ensure the optimality of $\bx^0$, the inverse problem needs to find the constraint parameters $\bA$ and $\bb$ such that $\bx^0$ is on the boundary of the imputed feasible region and can be optimal with respect to $\bc$. Hence, constraint \eqref{eq:TMIO-StrongDual} is necessary to ensure that strong duality holds for the preferred solution.}  In contrast to the single-point $\IO$ formulation, the primal feasibility constraint~\eqref{eq:TMIO-PrimalFeas} is now a set of $K$ constraints that ensures the feasibility of all observations for $\FO$. The formulation of $\MIO$, similar to that of $\IO$, is bilinear and hence, is non-convex in general. Analogous to $\IO$, we first show that the $\MIO$ formulation is feasible. 

\begin{proposition}\label{prop:MIOfeas}
\add{The feasible region of $\MIO$ is non-empty.} 
\end{proposition} 
\proof{Proof.}  
Let $\bw = \bzero$, $\by = || \bc||\bone$, $\bb = (\bc'\bx^0)/\|\bc\|$, and $\ba_i= \bc/\|\bc\|$ $\forall i \in \cI$. 
The resulting solution $(\bA, \bb, \by, \bw)$ satisfies constraints~\eqref{eq:TMIO-StrongDual}--\eqref{eq:MIO-Dualfeas2}. To show that the primal feasibility constraints~\eqref{eq:TMIO-PrimalFeas} also hold, note that if $\exists k \in \cK$ such that $\bA \bx^k < \bb$, then by substituting the values of $\bA$ and $\bb$, we have $(\bc' \bx^k)/\|\bc\| < (\bc'\bx^0)/\|\bc\|$, or equivalently, $\bc' \bx^k < \bc' \bx^0$, which is a contradiction to the definition of $\bx^0$ (Definition~\ref{def:x0}). Therefore, the  solution $(\bA, \bb, \by, \bw)$ is feasible for $\MIO$.
\Halmos \endproof

Constraint parameters $\bA$ and $\bb$ described in Proposition~\ref{prop:MIOfeas} represent the half-space $\cC=\{ \bx \in \mathbb{R}^n |~ \bc'\bx \geq \bc'\bx^0 \}$ whose identifying hyperplane is orthogonal to the cost vector $\bc$ and passes through the preferred solution $\bx^0$ (as an example, see Figure~\ref{fig:noPriorB:a} in Section~\ref{sec:numericalexample}). Therefore, the set $\cC$ includes all observations $\bx^k, k\in\cK$ (\ie, $\cC$ is a valid feasible set) and makes $\bx^0$ optimal for the forward problem. Hence, $\cC$ is a feasible set for the $\FO$ problem that is {\it imputed} from a solution of $\MIO$. In Definition~\ref{def:imputedSet}, we generalize this concept for all valid feasible sets that are derived from $\MIO$ solutions. 

\begin{definition} \label{def:imputedSet}
A polyhedron $\cX = \{\bx\in \mathbb{R}^n |~ \bD \bx \geq \bd \}$  is called an \underline{imputed} feasible set if there exists a feasible solution $(\bA, \bb, \by, \bw)$ of $\MIO$ such that $\cX = \{\bx\in \mathbb{R}^n |~ \bA \bx \geq \bb, \bG \bx \geq \bh \}$.
\end{definition}
An imputed feasible set $\cX = \{\bx\in \mathbb{R}^n \,|\, \bD \bx \geq \bd \}$ may be represented by infinitely many sets of constraints. For example, any scalar multiplication of the inequality or any other (perhaps linearly-dependent) reformulation  will represent the same set $\cX$. The $\MIO$ formulation finds one such set of constraints to characterize $\cX$ while satisfying the normalization constraint~\eqref{eq:MIOnorm}. When referring to an imputed feasible set, we consider the set $\cX$ and not the exact constraint parameters that define it. Note that any imputed feasible set is always a feasible region for $\FO$ that makes $\bx^0$ optimal because it is inferred by a solution of $\MIO$, and conversely, any feasible region of $\FO$ that satisfies the known constraints and makes $\bx^0$ optimal is an imputed feasible set of $\MIO$. We formalize this property in Proposition~\ref{prop:valid}. 

\begin{proposition}\label{prop:valid}  
The system $(\bA, \bb, \by, \bw)$ is feasible for $\TMIO$ if and only if the polyhedron $\cX = \{\bx \in \mathbb{R}^n|~ \bA \bx \geq \bb, \bG \bx \geq \bh\}$ is a valid feasible set that makes $\bx^0$ optimal for the forward problem.  
\end{proposition}
\proof{Proof.} 
Assume that $(\bA, \bb, \by, \bw)$ is a solution to $\TMIO$. By constraint~\eqref{eq:TMIO-PrimalFeas}, ${\bA \bx^k \geq \bb, \forall k\in \cK}$, and by Definition~\ref{assum:G}, we have $\bG \bx^k \geq \bh, \forall k\in \cK$. Hence, $\cX$ is a valid feasible set. Constraint~\eqref{eq:TMIO-StrongDual} ensures that strong duality holds for $\bx^0$, and hence, $\bx^0$ must be optimal for $\FO$.

\noindent Now let $\cX = \{\bx \in \mathbb{R}^n|~ \bA \bx \geq \bb, \bG \bx \geq \bh\}$ be a valid feasible set that makes $\bx^0$ optimal for $\FO$. Without loss of generality, we can assume that constraint~\eqref{eq:MIOnorm} holds since we can always normalize $\bA$ and $\bb$ so that $\|\ba_i\|=1$. The primal feasibility constraint~\eqref{eq:TMIO-PrimalFeas} is always met by definition of $\cX$. Since $\bx^0$ is optimal for $\FO$, we have $\underset{\bx \in \cX}{\min}\{\bc'\bx\}>-\infty$, and therefore, the dual of $\FO$ exists and is feasible, and strong duality holds. Hence, all constraints~(\ref{eq:TMIO-PrimalFeas}-\ref{eq:MIO-Dualfeas2}) are satisfied, which implies that there must exist $\by$ and $\bw$ such that $(\bA, \bb, \by, \bw)$ is feasible for $\MIO$.
\Halmos \endproof

Proposition~\ref{prop:valid} characterizes the properties of all solutions to $\MIO$ and ensures that $\bx^0$ is optimal for $\FO$. Although Proposition~\ref{prop:valid} and the $\MIO$ formulation explicitly consider the optimality of only $\bx^0$, we show in  Remark~\ref{rem:StrongDual} that any other observation with the same objective function value as $\bx^0$ is also optimal for the forward problem. 
\begin{remark}\label{rem:StrongDual} 
If $\cX$ is an imputed feasible set, any $\tilde{\bx} \in \cX$ such that $\tilde{\bx} \in \argmin_{\bx^k, \forall k\in \cK}\{\bc'\bx^k\}$ is an optimal solution of $\FO$.
\end{remark}
\proof{Proof.}
Let $\bx^0$ be the preferred solution. Assume $\exists \, \tilde{\bx} \in \cX$ such that $\tilde{\bx} \in \argmin_{\bx^k, \forall k\in \cK}\{\bc'\bx^k\}$ and ${\tilde{\bx} \ne \bx^0}$. By constraint~\eqref{eq:TMIO-PrimalFeas}, we know that $\tilde{\bx}$ is feasible for $\FO$. If $\tilde{\bx}$ is not optimal for $\FO$, then $\bc' \bx^0 < \bc' \tilde{\bx}$ which is a contradiction to $\tilde{\bx} \in \argmin_{\bx^k, \forall k\in \cK}\{\bc'\bx^k\}$. Hence, $\tilde{\bx}$ must be an optimal solution to $\FO$.
\Halmos \endproof

\add{In this section, we proposed inverse optimization models that can impute the feasible region of a forward problem based on a set of feasible observations, and we discussed the general properties of the solutions. Considering that the proposed models are nonlinear, we next focus on additional properties of the solutions and propose a tractable reformulation that can be used to solve the $\MIO$ problem.}

\subsection{Tractable Reformulation}

\add{The proposed $\MIO$ formulation includes a set of bilinear constraints which makes the formulation non-linear (\ie, constraints~\eqref{eq:TMIO-StrongDual} and~\eqref{eq:TMIO-DualFeas}) and therefore, intractable to solve. In this section, we outline specific properties of the solution space of the $\MIO$ model that would allow us to develop a tractable reformulation of this model. To this end, we first characterize the range of possible imputed feasible sets to $\MIO$ and then prove that by considering specific known constraints, strong duality and dual feasibility can be guaranteed without explicitly incorporating the corresponding nonlinear constraints in the formulation. We finally formalize this idea theoretically and discuss how this reformulation can be used to find solutions to $\MIO$.}

We first show that we can find the smallest and largest possible imputed feasible sets of $\MIO$ solely based on the given observations. An imputed feasible set $\cX$ of $\MIO$ is a valid feasible set, according to Proposition~\ref{prop:valid}. By Remark~\ref{rem:validSet}, any superset of $\cX$ will also be a valid feasible set. In particular, $\cX$ is always a subset of the half-space $\cC = \{ \bx \in \mathbb{R}^n |~ \bc'\bx \geq \bc'\bx^0 \}$ and a superset of the convex hull of the observations, denoted by $\cH$. This property is shown in Lemma~\ref{lem:subset} and plays a fundamental role in reformulating the $\MIO$ model later in Theorem~\ref{thm:MIP_general}. For brevity of notations, we use $\cH$ and $\cC$ as defined above throughout the rest of this paper.

\begin{lemma}\label{lem:subset}
If $\cX$ is an imputed feasible set of $\MIO$, then $\cH \subseteq \cX \subseteq \cC$. 
\end{lemma}
\proof{Proof.}
\noindent ($\cH\subseteq \cX$): Assume $\cH\not\subseteq \cX$ and $\exists\, \bar{\bx} \in \cH,$ $\bar{\bx} \not \in \cX$. By definition of $\cH$, $\exists \, \lambda_k \geq 0, \, \forall k \in \cK$ such that $\bar{\bx} = \sum_{k \in \cK} \, \lambda_k \bx^k$ and $ \sum_{k \in \cK} \, \lambda_k = 1$. This is a contradiction because $\cX$ is a polyhedron that is a valid feasible set. Therefore, it contains all observations $\bx^k$ and any convex combination of them, including $\bar\bx$. Hence, $\cH\subseteq \cX$. \\
($\cX \subseteq \cC$): Similarly, assume $\cX \not \subseteq \cC$ and $\exists \, \bar{\bx} \in \cX,$ $\bar{\bx} \not \in \cC$. Since $\bar{\bx} \not \in \cC$, we have $\bc'\bar{\bx} < \bc' \bx^0$ (by definition). Therefore, $\bar{\bx}$ has a better objective value than $\bx^0$, which is a contradiction to $\cX$ being an imputed feasible set because $\cX$ must make $\bx^0$ optimal for $\FO$. Hence, $\cX \subseteq \cC$. 
\Halmos \endproof

As Lemma~\ref{lem:subset} illustrates, any imputed feasible set must be a subset of the half-space $\cC = \{ \bx \in \mathbb{R}^n |~ \bc'\bx \geq \bc'\bx^0 \}$.  
\add{The intuition behind this idea is as follows. The identifying hyperplane of $\cC$ (\ie, $\bc'\bx = \bc'\bx^0 $) passes through the preferred solution $\bx^0$ and is orthogonal to the known cost vector $\bc$. The inclusion of this half-space ensures that $\bx^0$ is on the boundary of the imputed feasible region and is always candidate optimal. In other words, for a valid feasible set $\cU \not \subseteq \cC$, there will always exist other feasible solutions that have a better objective function value than $\bx^0$ in the forward problem. Such a set $\cU$ cannot be an imputed feasible set of $\MIO$ since it does not make $\bx^0$ a candidate for optimality. Hence, any imputed feasible set must be a subset of the half-space $\cC$. For a visual representation of the half-space $\cC$, see Figure~\ref{fig:noPriorB:a} in Section~\ref{sec:numericalexample}.} 
Using this property, we can reduce the solution space of $\MIO$ from $\mathbb{R}^n$ to the half-space $\cC$. 

We can further restrict the solution space of $\MIO$ by noting that the set of known constraint $\cG = \{ \bx \in \mathbb{R}^n \mid \bG \bx \geq \bh \}$ also has to be met for any $\MIO$ solution. Therefore, the solution space of $\MIO$ is always a subset of $\cS = \cC \cap \cG$. Proposition~\ref{prop:CG} implies that this solution space is the largest imputed feasible set of $\MIO$ and any $\MIO$ solution is a subset of the space $\cS$.

\begin{proposition}\label{prop:CG} 
\add{Let $\cS =\cC \cap \cG =  \{ \bx \in \mathbb{R}^n |~ \bc'\bx \geq \bc'\bx^0, \, \bG \bx \geq \bh\}$, then  
\begin{enumerate}   \setlength\itemsep{0em}
    \item [{\normalfont(}i\,{\normalfont)}] $\cS$ is an imputed feasible set of $\MIO$, 
    \item [{\normalfont(}ii\,{\normalfont)}] for any other imputed feasible set $\cX$, we have $\cX \subseteq \cS$,
    \item [{\normalfont(}iii\,{\normalfont)}] for any valid feasible set $\cU$, the set $\cU \cap \cS$ is an imputed feasible set. 
\end{enumerate}
}
\end{proposition}

\proof{Proof.}
\add{
({\it i\,}) The set $\cS$ is a valid feasible set since both $\cC$ and $\cG$ are valid feasible sets as shown in Lemma~\ref{lem:subset} and Assumption~\ref{assum:G}. The set $\cS$  also makes $\bx^0$ optimal for $\FO$ because it includes the half-space $\cC$. Hence, by Proposition~\ref{prop:valid}, $\cS$ is an imputed feasible set of $\MIO$. \\
({\it ii\,})  For any imputed feasible set $\cX$, it is obvious that $\cX \subseteq \cC$ (by Lemma~\ref{lem:subset}) and $\cX \subseteq \cG$ (by definition), and hence $\cX \subseteq  \cS$. \\
({\it iii\,}) Since both $\cU$ and $\cC$ are valid feasible sets of $\FO$, the set $\cU \cap \cS$ is also a valid feasible set, and hence, primal feasibility holds. Strong duality also holds for $\bx^0 \in \cU \cap \cS$ because the half-space $\cC$ is considered as part of $\cS$ with $\bc'\bx^0$ as the optimal value of the $\FO$. Hence, dual feasibility also holds, and $\cU \cap \cS$ is an imputed set of $\MIO$.
}
\Halmos \endproof

Without loss of generality, in the rest of this paper, we assume that $\cC=\{ \bx \in \mathbb{R}^n |~ \bc'\bx \geq \bc'\bx^0 \}$ is added as the {\it first} unknown constraint in the formulation, that is, $\bg_1 = \bc, \, h_1 = \bc' \bx^0$. Let the set $\cS =\cC \cap \cG =  \{ \bx \in \mathbb{R}^n |~ \bc'\bx \geq \bc'\bx^0, \, \bG \bx \geq \bh\}$ denote the new set of ``known constraints'' hereinafter. With this assumption, Remark~\ref{rem:cs} formally points out that the largest possible feasible set that can be imputed by a solution of $\MIO$ is the set $\cS$ itself, as defined in Proposition~\ref{prop:CG}.

\begin{remark}\label{rem:cs}
The set $\cS$ is the largest possible imputed feasible set of $\MIO$.
\end{remark}

Considering the set $\cS$ as the set of known constraints, we can guarantee that the strong duality and dual feasibility constraints~\eqref{eq:TMIO-StrongDual} and~\eqref{eq:TMIO-DualFeas} hold without explicitly including them in the model. Based on these properties, Theorem~\ref{thm:MIP_general} shows an equivalent reformulation of the $\MIO$ problem when the half-space $\cC$ is considered as a known constraint.

\begin{theorem}\label{thm:MIP_general} 
Solving $\MIO$ is equivalent to solving the following problem when \add{$\cS =\cC \cap \cG =  \{ \bx \in \mathbb{R}^n |~ \bc'\bx \geq \bc'\bx^0, \, \bG \bx \geq \bh\}$ is the set of known constraints} of $\FO$.
\begin{subequations}  \label{eq:MIP_general} 
\begin{align}
    {\mathbf \eMIO:~}\underset{\bA, \bb}
    {\normalfont{\text{minimize}}} \quad &  \mf(\bA, \bb; \bP) \\    
    {\normalfont{\text{subject to}}} \quad & \ba_{i}' \, \bx^k \geq b_i, \qquad \forall i \in \cI, 
    \quad  k \in \cK  
    \label{eq:mip1}\\
     & ||\ba_i|| = 1 \label{eq:MIPnorm},\\ 
    & {\color{myGreen} \bA \in \mathbb{R}^{m_1 \times n}, \quad \bb \in \mathbb{R}^{m_1}.} \label{eq:MIPboxConst} 
\end{align} 
\end{subequations} 
\end{theorem}
\proof{Proof.} 
({\it i}\,) If $(\bA, \bb, \by, \bw)$ is a solution of $\MIO$, then  $(\bA, \bb)$ is a solution to the $\eMIO$ formulation since~\eqref{eq:mip1}--\eqref{eq:MIPboxConst} are also constraints of  $\MIO$. 
({\it ii}\,) Conversely, for the pair $(\bA,\bb)$ that is a solution to $\eMIO$, let $\bw = (1, 0, \dots, 0)$, $\by = (0, 0, \dots, 0)$. The solution $(\bA, \bb, \by, \bw)$ is feasible for $\MIO$ since by Proposition~\ref{prop:CG}, the strong duality constraint \eqref{eq:TMIO-StrongDual} and the dual feasibility constraint \eqref{eq:TMIO-DualFeas} hold through the inclusion of the half-space $\cC$ as a known constraint in $\cS$. Therefore, by ({\it i}\,) and ({\it ii}\,), solving $\eMIO$ is equivalent to solving $\MIO$. 
\Halmos \endproof

Theorem~\ref{thm:MIP_general} shows that by considering the half-space as one of the known constraints, instead of solving the bilinear $\MIO$ problem, we can solve a simpler problem that does not explicitly include the strong duality and dual feasibility constraints and hence, does not have any bilinear terms. Note that there are multiple ways to re-write constraint~\eqref{eq:MIPnorm} based on the particular application and the desired properties of the resulting model. Depending on the type of normalization constraint used in~\eqref{eq:MIPnorm}, the complexity of the corresponding $\eMIO$ formulation would differ. For example, popular norms such as $L_1$ or $L_2$ would yield linearly- or quadratically-constrained problems, respectively.  

Remark~\ref{rem:eMIO} highlights that any valid feasible set for $\FO$ can be a solution to $\eMIO$. This property is intuitive since the $\eMIO$ formulation only includes the primal feasibility constraints for all observations. Hence, the feasible region of $\eMIO$ reduces to the set of valid feasible sets of $\FO$. Therefore, by the inclusion of $\cC$ in the known constraints, the complexity of the problem reduces to only finding valid feasible sets of $\FO$ through the $\eMIO$ formulation. 
\begin{remark}\label{rem:eMIO}
Any valid feasible set $\cX$ of $\FO$ is an imputed feasible set to $\eMIO$. 
\end{remark}
\proof{Proof.}
The set $\cX$ is a valid feasible set, and the set of known constraints in $\eMIO$ is $\cS$. Therefore, by Proposition~\ref{prop:CG}, the set $\cX \cap \cS$ is an imputed feasible set to $\eMIO$.  
\Halmos \endproof

We finally note that solving $\eMIO$ provides constraint parameters $\bA$ and $\bb$ such that $\bA\bx \geq \bb$ along with the set of known constraints $\cS$ shape the imputed feasible region of $\FO$. In other words, any solution to $\eMIO$ can identify an imputed feasible set of $\MIO$ by first finding the constraint parameters $\bA$ and $\bb$ and then finding the intersection of these constraints with the known constraints. Remark~\ref{cor:eMIO} formalizes this concept.  

\begin{remark}\label{cor:eMIO}
An imputed feasible set of $\MIO$ can be derived as $\{\bx \in \mathbb{R}^n |~ \bA \bx \geq \bb, \bx \in \cS\}$ for any solution $(\bA, \bb)$ to $\eMIO$. 
\end{remark}

In this section, we showed that the solution to the $\MIO$ (or $\IO$) formulation can be found by first solving the $\eMIO$ formulation and then deriving the corresponding imputed feasible set through adding the known constraints. The complexity of $\eMIO$ depends on the type of norm in constraint~\eqref{eq:MIPnorm} and the complexity of the loss function. The $\eMIO$ problem will be a linearly-constrained model if a linear norm is used. In the next section, we focus on different types of loss functions and provide specific examples of measures to induce the characteristics of the feasible region.

\section{Loss Functions} \label{sec:measures}
The $\MIO$ formulation minimizes an objective function $\mf(\bA, \bb, \bP)$ which affects the optimal solution $(\bA, \bb, \by, \bw)$ and hence, drives the desirable properties of the imputed feasible set of $\FO$. This imputed feasible set for the forward problem may take various shapes and forms based on the given parameter set ($\bP$) and the objective function  ($\mf$). In this chapter, we introduce several loss functions that can be used based on the available information on the constraints. 
Note that all the models introduced in our framework include a generalized loss function $\mf$, and this function can be tailored by the user to induce properties for the application domain at hand. As shown in Section~\ref{sec:Methodology}, the solutions to both the $\IO$ and $\MIO$ formulations can be found by solving the $\eMIO$ model. Therefore, without loss of generality, we use the $\eMIO$ formulation to develop the theoretical properties of models with different loss functions in this section.

In the literature of inverse optimization, a \emph{prior belief} on the constraint parameter is defined as reasonable or desired values for the constraint parameters, and the inverse problem often attempts to minimize the distance of the parameters from such belief. In our framework, we do not necessarily require the user to provide any such prior belief on the parameters. Therefore, for any unknown constraint in $\FO$, we consider two cases:  ({\it i}\,)  a prior belief on the constraint is available, and  ({\it ii}\,) no such prior information exists. For case ({\it i}\,), which has been considered in the literature, we discuss a specialization of our general loss function that allows the user to minimally perturb these prior beliefs. In case ({\it ii}\,), where no information on the constraint is assumed, we show that it is possible to find a large variety of imputed feasible sets for $\MIO$. We introduce different loss functions that aim to find the appropriate constraints when no prior belief on the constraint parameters is available.

In the rest of this section, we first discuss the theoretical properties of imputed feasible sets of $\eMIO$ when a prior belief on the constraints is available. Next, we present and discuss different loss functions that can be employed in the absence of a prior belief.

\subsection{Prior Belief on Constraints Available} \label{sec:priorB}
When a prior belief on the constraint parameters is available, the objective of the inverse problem is often to minimize some measure of distance (\eg, norm) of the imputed constraint parameters from that prior belief. In this section, we study the use of prior belief as a loss function in the objective of our $\eMIO$ model. We refer to this loss function as the {\it \priorB}. 
Let the assumed prior belief on the constraint parameters, denoted as $\hatA$ and $\hatb$, be given as the input parameter $\bP$. For ease of notations, let $\barA = [\bA \,\,\,\bb]$ be the matrix that appends the column $\bb$ to the $\bA$ matrix and $\barhatA = [\hatA \,\,\, \hatb]$ be the corresponding prior belief. We define the \priorB~as the loss function that captures the distance of $\barA$ from the prior belief $\barhatA$ according to some norm $||\cdot ||$ as follows: 
\[\mf(\bA, \bb; \bP) = \mf(\bA, \bb, \barhatA) = \sum_{i \in \cI}\omega_i\, ||\barA_i - \barhatA_i ||. \qquad \tag{\priorB}\]
Parameter $\omega_i$ is the objective weight capturing the relative importance of constraint $i$, and $\barA_i$ and $\barhatA_i$ are the $i^\text{th}$ rows of matrices $\barA$ and $\barhatA$, respectively. 
Proposition~\ref{prop:prior2} shows that the $\eMIO$ model with the \priorB~can be decomposed into solving a series of smaller problems for each of the $m_1$ unknown constraints. 

\begin{proposition}\label{prop:prior2}
The optimal solution of $\eMIO$ with the \priorB~can be found by solving the following problem $m_1$ times for each $i \in \cI = \{1,\dots,m_1 \}$. 
\begin{subequations} \label{eq:prior}
\begin{align}
    \underset{\ba_i, b_i }{\normalfont\text{minimize}} \quad & ||\barA_i - \barhatA_i || \\
{\normalfont\text{subject to}} \quad    & \ba_{i}' \, \bx^k \geq b_i, \quad \forall k \in \cK
    \label{eq:prior1i}\\
    & ||\ba_i|| = 1. 
\end{align}
\end{subequations}
\end{proposition}
\proof{Proof.} The $\eMIO$ problem with the \priorB~is separable for each constraint $i$, which means problem~\eqref{eq:prior} can be solved $m_1$ times to recover each $\ba_i$ and $b_i$ independently.
\Halmos \endproof

If the prior belief is not a valid feasible set, then at least one of the observations $\bx^k, \, k \in \cK$ is positioned outside of the prior belief. Therefore, $\barhatA$ needs to be minimally perturbed to generate a valid feasible set. This is a prevalent occurrence in practice since although a set of {\it a priori} constraints might be available, in reality, these constraints might be too tight to hold for all observations.  
If the set identified by the prior belief $\barhatA$ is a valid feasible set, \ie, $\hatA \bx^k \geq \hatb, \, \forall k \in \cK$, then Proposition~\ref{prop:priorValidSet} shows that there is a closed-form solution to $\eMIO$. 
\begin{proposition}\label{prop:priorValidSet}
If $\cX = \{\bx\in \mathbb{R}^n \, |~  \hatA \bx \geq \hatb \}$ is a valid feasible set, then $\bA = \hatA$ and $\bb = \hatb$ is an optimal solution to $\eMIO$ under the \priorB. 
\end{proposition}
\proof{Proof.}
By assumption, $\cX$ is a valid feasible set and hence by Remark~\ref{rem:eMIO}, an imputed feasible set to $\eMIO$. Therefore, $\barA = \barhatA$ is a feasible solution to $\eMIO$ with $\mf(\bA, \bb; \bP) = 0$ under the \priorB. Hence, $\bA = \hatA$, $\bb = \hatb$ is an optimal solution to $\eMIO$. 
\Halmos \endproof

The \priorB, which is most often used in the literature, heavily relies on both the availability and the quality of the prior belief. In particular, if the quality of the prior belief is poor, it enforces the inverse optimization to fit the imputed feasible set to this poor-quality prior belief. In what follows, we propose and discuss other loss functions that can be employed if no quality prior belief is available for the constraint parameters.

\subsection{No Prior Belief on Constraints} \label{noPriorB}
In this section, instead of relying on a prior belief, we propose different loss functions that can incorporate other data (\eg, observations) to find the solution of $\eMIO$. We start with a simple constraint satisfaction model, which is sometimes used in the literature of inverse optimization. We then propose three new loss functions that each result in different properties for the imputed feasible set of $\eMIO$. We also consider combining the loss functions to further refine the shape of the imputed feasible region. 

\subsubsection{\zeroObj}
If no preference and no information about the feasible region is given, \ie, there is no data provided to be used to derive the shape of the imputed feasible set ($\bP=[\,\,]$), then the $\eMIO$ reduces to a feasibility problem by setting the objective function as zero, \ie, \[\mf(\bA, \bb; \bP) = 0.  \tag{\zeroObj} \] 
We refer to this loss function as the {\it \zeroObj}.  
\begin{proposition}\label{prop:feassol}
A closed-form optimal solution for $\eMIO$ with the \zeroObj~is
\begin{align}
  \ba_{i} = \frac{c_i}{||\bc||}, 
  \qquad  
  b_i = \frac{\bc'\bx^0 }{||\bc||}, \qquad \forall i \in \cI. 
\end{align}
\end{proposition}
\noindent As shown in Section~\ref{sec:Methodology}, the solution above is feasible for $\eMIO$ and is hence, optimal under the \zeroObj. Intuitively, any feasible solution to $\eMIO$ is an optimal solution in this case. This property is highlighted in Remark~\ref{prop-infinitefeas}. 

\begin{remark}\label{prop-infinitefeas}
The $\eMIO$ formulation with the \zeroObj~has an infinite number of optimal solutions. 
\end{remark}
\proof{Proof.} 
The convex hull $\cH$ of the observations is a valid feasible set by definition, and hence, it is an imputed feasible set of $\eMIO$ under the \zeroObj. By Remark~\ref{rem:validSet}, any set $\cX$ that $\cH \subseteq \cX$ is also a valid feasible set, and by Proposition~\ref{prop:CG}, $\cX\cap \cC$ is an imputed feasible set for $\MIO$ (and hence, for $\eMIO$). Therefore, $\eMIO$ has infinitely many imputed feasible sets, and accordingly, infinitely many optimal solutions.
\Halmos \endproof

In practice, the \zeroObj~may not be the loss function of choice if there exist some properties that are preferred for the feasible set of $\FO$. In the rest of this section, we introduce three other loss functions that can inform the shape of the imputed feasible set using the observations and discuss their properties. 

\subsubsection{\minDist} 
The \emph{\minDist}~finds a feasible region that has the smallest total distance from all of the observations. Here, the given parameter $\bP$ is the matrix that includes all observations, $\bP=[\bx^1,\dots,\bx^k]$. This  loss function minimizes the sum of the distances of each observation from all constraints. Let $d_{ik}$ denote the distance of each observation $\bx^k, \, k\in \cK$ from the identifying hyperplane of the $i^{\text{th}}$ constraint. The {\it \minDist}~is defined as 
\[\mf(\bA, \bb, \bP) = \mf(\bA, \bb, [\bx^1,\dots,\bx^k])=\sum_{k=1}^{K} \sum_{i=1}^{m_1} d_{ik}, \quad \tag{\minDist}\]
where the distance $d_{ik}$ can be calculated using any distance metric, for example, the \emph{Euclidean distance}, or the \emph{slack distance} defined as $d_{ik}~=~\ba_{i} \bx^{k}- b_{i}$. Similar to the \priorB, this loss function is separable for each constraint, and hence, the resulting $\eMIO$ model can be decomposed and solved for each constraint independently, as shown in Proposition \ref{prop:mindist}.  

\begin{proposition}\label{prop:mindist}
The optimal solution of $\eMIO$  with the \minDist~can be found by solving  the following problem $m_1$ times, for each $i\in\cI$: 
 \begin{subequations}
 \begin{align}
   \underset{\ba_i, b_i}  {\normalfont\text{minimize}} \quad &\sum_{k=1}^{K}
   d_{ik} \\ 
 {\normalfont\text{subject to}} \quad & \ba_{i}' \, \bx^k \geq b_i, \quad \forall k \in \cK
    \label{eq:mindist:i}\\
     & ||\ba_i|| = 1.   
 \end{align}
 \end{subequations}
\end{proposition}

\subsubsection{\fairness} 
This loss function aims to find a feasible set {such that all of its constraints are} equally close to all observations and hence, is ``fair''. Using the same notations as those in the \priorB, we calculate  $d_{ik}$ as the distance of each observation $k$ from the identifying hyperplane of each constraint $i$. We then calculate the \add{total} distance for all observations, $d_k=\sum_{i=1}^{m_1} d_{ik}$. The {\it \fairness}~is 
\[\mf(\bA, \bb, \bP)  = \mf(\bA, \bb; [\bx^1,\dots,\bx^K])=\sum_{k\in\cK}{(d_k-\sum_{k\in \cK}d_k/K)}. \quad \tag{\fairness} \]

This measure minimizes the deviation of the total distances for all observations and ensures that all observations have roughly the same total distance from all constraints. The \fairness~avoids cases were the constraints are all on one side of the observations and far away from others, and hence, it typically results in imputed feasible sets that are more confined compared to the \minDist. 

\subsubsection{\minmin} 
The {\it \minmin}~tries to find the constraint parameters such that the minimum distance of each observation from all of the constraints is minimized. In other words, it tries to ensure that each observation is close to at least one constraint, if possible (\ie, if the observation is not an interior point). Again, let $d_{ik}$ be the distance of observation $k$ from the identifying hyperplane of constraint $i$. The \minmin~is defined as
\[\mf(\bA, \bb, \bP)  = \mf(\bA, \bb; [\bx^1,\dots,\bx^k])=\sum_{k\in\cK}{ \min_{i\in\cI} d_{ik}}. \quad \tag{\minmin}\]
Minimizing the \minmin~can be written as $\, \min \sum_{k\in\cK}{ \min_{i\in\cI} d_{ik}}$, and this min-min objective can be reformulated using auxiliary binary variables. The resulting $\eMIO$ formulation under the \minmin~is as follows:
\begin{subequations}
\begin{align}
 \text{minimize} \quad  & \sum_{k \in \cK} m_k \label{eq:minmin_obj}\\
 \text{subject to} \quad  & d_{ik} =  \sum_{j \in \cJ}a_{ij} x^k_{i} - b_i  \quad \forall i\in \cI, \, k \in \cK \label{eq:minmin_const_first}\\
 &  m_k \geq d_{ik} - M \gamma_{ik} ,   \quad \forall i\in \cI, \, k \in \cK\\ 
 &  \sum_{i} \gamma_{ik} = |\cI| - 1, \quad \forall k \in \cK \\
 & \gamma_{ik} \in \{0,1\}, \quad \forall i \in \cI, k \in \cK \\
 & m_k \geq 0, \quad \forall k \in \cK  \label{eq:minmin_const_last}\\
 &\eqref{eq:mip1}-\eqref{eq:MIPboxConst}.
\end{align}
\end{subequations}

Note that the resulting model is a mixed-integer linear program if a linear norm is used as constraint~\eqref{eq:MIPnorm} of $\eMIO$. 

\subsubsection{Combined Loss Functions} 
In the literature, inverse optimization formulations tend to produce multiple optimal solutions and return one of them arbitrarily as the optimal solution. Our inverse optimization formulations often demonstrate this property as well, even when the previously-mentioned loss functions are imposed. Given that the $\eMIO$ formulation is tractable, we can utilize the multi-optimum property of inverse optimization to further calibrate the shape and characteristics of the imputed feasible set of $\FO$ by combining multiple loss functions. 

Typical approaches for combining different objective functions include using multi-objective optimization and using sequential objectives. The former is trivial to implement but introduces challenges such as deciding on the weights of each of the multiple objectives and is still prone to generating multiple optimal solutions that do not necessarily reflect the desired characteristics. The latter approach (also referred to as {\it secondary objective}) is our suggested method since each iteration narrows down the solution space to further fine-tune the solution to the specific characteristics of interest. We note that this approach does not require significant additional computational burden given that the $\eMIO$ formulation is linearly-constrained for some popular norms (\eg, $L_1$) and multiple instances can be solved sequentially.

In the secondary objective approach, the model is solved for a loss function of choice, say $\mf_1$, where the optimal value of $\mf_1^*$ is achieved. Then, to select those imputed feasible sets whose corresponding solutions generate the same optimal value of $\mf_1^*$ but possess other desired properties as well, the $\eMIO$ is solved again with a new loss function $\mf_2$ and an additional  constraint of $\mf_1(\bA, \bb, \bP) = \mf_1^*$. This process can be repeated for as many loss functions as desired. 

In the next section, we test our approach on two numerical examples and compare the results for different loss functions. We then provide an example of combining the loss functions, by choosing the \fairness~as the primary objective $\mf_1$ and the \minDist~as the secondary objective $\mf_2$. In this case, we find solutions with the best \fairness~value that are also closer to all observations with regards to the \minDist.

\section{Numerical Results} \label{sec:numericalexample}
In this section, we test our methodology on two illustrative two-dimensional (2D) numerical examples \add{and a larger-scale diet recommendation case study}. For the ease of visualization, we use the 2D datasets ($n=2$) to graphically show the observations, the feasible region, and the objective vector. In the first example, we consider a small number of equidistant observations ($K=5$) that form a symmetric convex hull. For this example, the inverse solutions are easy to find by visual inspection, which allows us to understand the intuition behind the solutions generated under each loss function and compare their characteristics. The second numerical example considers a relatively larger set of observations ($K = 19$) that are randomly placed and their convex hull has an arbitrary shape. This example further elaborates on the insights from each of the introduced loss functions under non-trivial cases. In each of the two examples, we consider multiple known and unknown constraints in the $\FO$ problem. We solved the first example using both the $\MIO$ and $\eMIO$ formulations which confirmed the equivalence of the results for the two models. The second example, however, was only solved using the $\eMIO$ model because the commercial solver we used was not able to solve the larger non-linear model to optimality. \add{Finally, we apply the $\eMIO$ formulation to a much larger example of a diet recommendation application. In this case study, we consider $K=100$ observations of a dieter's daily food intake from a set of $n=26$ food items. We consider a set of known nutrient constraints and impute multiple implicit constraints of the dieter. We compare the palatability of the resulting diet recommendations with and without the imputed constraints. }

In our numerical results, we use the $L_2$~norm in the \priorB~since it is a popular norm used in the literature. For all other loss functions, for simplicity, we use the linear slack distance (\ie, $d_{ik}=\ba_i'\, {\bx^k} - b_i$) to calculate the distance of a given point $\bx^k$ to the identifying hyperplane of the $i^{\text{th}}$ constraint (\ie, $\ba_i' \, \bx = b_i$). We note that there exist other linear distance metrics (\eg, $L_{\infty}$~norm) that can be used, but we find the slack distance to be more illustrative in a two-dimensional setting.
For the normalization constraint~\eqref{eq:MIOnorm}, we use $|\sum_{j\in \cJ}a_{ij}|=1$ as a proxy for the $L_1$~norm (\ie, $\sum_{j\in \cJ}|a_{ij}|=1$). We chose this normalization method instead of the $L_1$~norm to reduce the number of auxiliary binary variables to only $2n$ (as opposed to $2n\,(m_1+m_2)$) when reformulating it as a linear mixed-integer model. 

\subsection{Numerical Case I}\label{sec:NumCase1}
In the first numerical case, we have 5 observations, as listed in Table~\ref{tab:case1}. There are two known constraints with the first one being the half-space  $\cC$, as discussed in Section~\ref{sec:Methodology}.  
The $\MIO$ {and $\eMIO$} models were solved using the nonlinear solver {MINOS}~(\citeyear{saunders2003minos}) Version~5.51 and {CPLEX}~(\citeyear{cplex}) Version~12.9, respectively. Both models were formulated using AMPL~(\citeyear{ampl}) modeling language Version~20190529. While MINOS and other nonlinear solvers are sometimes capable of solving small-scale instances to optimality, they often fail to provide a global optimal solution in larger cases. For this numerical example, {MINOS} was able to solve the $\MIO$ model to optimality, and the $\MIO$ and $\eMIO$ solutions confirmed the same solutions in all instances. 

\begin{table}[htbp]
    \centering
 \begin{tabular}{l l}
    \toprule
    {\bf Description} & {\bf Value(s)} \\
    \midrule 
    Cost vector $\bc$ & $(-1, -1)$ \\
    \midrule
    Observations $\bx^0$; $\bx^k \qquad$  &  $\mathbf{(2,2)}$; $(1,1)$, $(1,2)$, $(2,1)$, $(1.5,1.5)$\\
    \midrule
    Known constraints   
                        & $0.5 x_1 + 0.5 x_2 \leq 2\quad$  (half-space $\cC$) \\ 
                        &  $x_1 + x_2 \geq 1 $ \\
    \midrule
    Unknown constraints & 4 constraints\\
    \bottomrule
\end{tabular}
    \caption{Numerical Case I}
    \label{tab:case1}
\end{table}

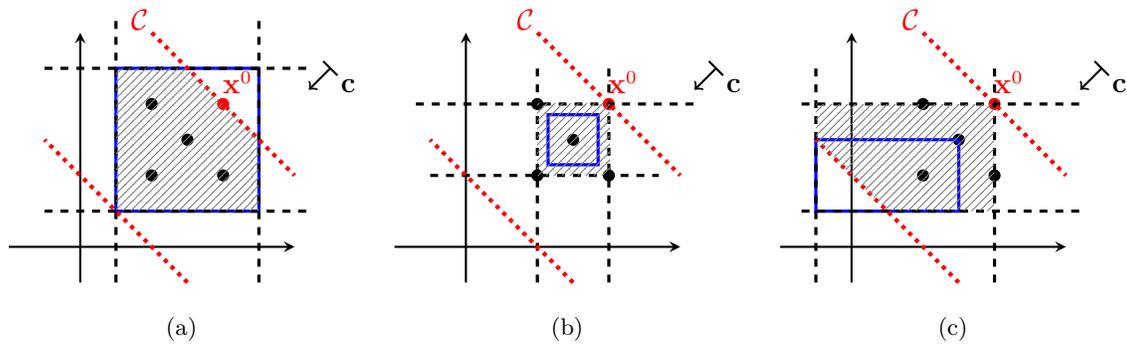
\begin{figure}[htbp] 
\centering
\subfigure[\label{fig:priorB:a}]{   
\begin{tikzpicture} [thick, scale=.95]
\draw[->, >=stealth] (-1,0) -- (3,0);
\draw[->, >=stealth] (0,-.5) -- (0,3);
\draw [fill, color=black] (1.5,1.5) circle [radius=0.07];
\draw [fill, color=black] (1,1) circle [radius=0.07];
\draw [fill, color=black] (2,1) circle [radius=0.07];
\draw [fill, color=black] (1,2) circle [radius=0.07];
\draw [fill, color=red] (2,2) circle [radius=0.07];
\node[anchor=south, color=red] at (2.2 , 2) (coord) {$\bx^0$};
\draw[<-] (3.2,2.2) coordinate -- (3.5,2.5) node[anchor= north west] {$\bc$};
\draw[-] (3.6, 2.4) -- (3.4, 2.6);
\draw[color=blue, very thick] (0.5,0.5) -- (0.5,2.5) -- (2.5,2.5) -- (2.5, 0.5) -- cycle;
\draw[dashed, very thick] (-.5, .5) -- (3.2, .5);
\draw[dashed, very thick] (-.5, 2.5) -- (3.2, 2.5);
\draw[dashed, very thick] (.5, -.5) -- (.5, 3.2);
\draw[dashed, very thick] (2.5, -.5) -- (2.5, 3.2);
\draw[dotted, ultra thick, red] (1, 3) -- (3, 1);
\node[anchor=south east, color=red] at (1.1,2.9) (coord) {$\cC$};
\draw[dotted,  ultra thick, red] (1.5,-0.5 ) -- (-0.5, 1.5);
\fill[pattern=north east lines, pattern color=gray] (.5, .5) -- (.5,2.5) -- (1.5, 2.5) -- (2.5, 1.5) -- (2.5,.5) -- cycle;
\end{tikzpicture}
}
\subfigure[\label{fig:priorB:b}]{   
\begin{tikzpicture} [thick, scale=.95]
\draw[->, >=stealth] (-1,0) -- (3,0);
\draw[->, >=stealth] (0,-.5) -- (0,3);
\draw [fill, color=black] (1.5,1.5) circle [radius=0.07];
\draw [fill, color=black] (1,1) circle [radius=0.07];
\draw [fill, color=black] (2,1) circle [radius=0.07];
\draw [fill, color=black] (1,2) circle [radius=0.07];
\draw [fill, color=red] (2,2) circle [radius=0.07];
\node[anchor=south, color=red] at (2.2 , 2) (coord) {$\bx^0$};
\draw[<-] (3.2,2.2) coordinate -- (3.5,2.5) node[anchor= north west] {$\bc$};
\draw[-] (3.6, 2.4) -- (3.4, 2.6);
\draw[color=blue, very thick] (1.15,1.15) -- (1.15, 1.85) -- (1.85, 1.85) -- (1.85, 1.15) -- (1.15, 1.15) -- cycle;
\draw[dashed, very thick] (-0.5, 1) -- (2.7, 1);
\draw[dashed, very thick] (-.7, 2) -- (3.2, 2);
\draw[dashed, very thick] (1, -.5) -- (1, 2.5);
\draw[dashed, very thick] (2, -.5) -- (2, 2.5);
\draw[dotted, ultra thick, red] (1, 3) -- (3, 1);
\node[anchor=south east, color=red] at (1.1,2.9) (coord) {$\cC$};
\draw[dotted, ultra thick, red] (1.5,-0.5 ) -- (-0.5, 1.5);
\fill[pattern=north east lines, pattern color=gray] (1,1)--(1,2)--(2,2)--(2,1)--cycle;%
\end{tikzpicture}
}
\subfigure[\label{fig:priorB:c}]{   
\begin{tikzpicture} [thick, scale=.95]
\draw[->, >=stealth] (-1,0) -- (3,0);
\draw[->, >=stealth] (0,-.5) -- (0,3);
\draw [fill, color=black] (1.5,1.5) circle [radius=0.07];
\draw [fill, color=black] (1,1) circle [radius=0.07];
\draw [fill, color=black] (2,1) circle [radius=0.07];
\draw [fill, color=black] (1,2) circle [radius=0.07];
\draw [fill, color=red] (2,2) circle [radius=0.07];
\node[anchor=south, color=red] at (2.2 , 2) (coord) {$\bx^0$};
\draw[<-] (3.2,2.2) coordinate -- (3.5,2.5) node[anchor= north west] {$\bc$};
\draw[-] (3.6, 2.4) -- (3.4, 2.6);
\draw[color=blue, very thick] (-0.5,0.5) -- (1.5,0.5) -- (1.5,1.5) -- (-.5, 1.5) -- (-.5, 0.5) -- cycle;
\draw[dashed, very thick] (-.7, .5) -- (3.2, .5);
\draw[dashed, very thick] (-.7, 2) -- (3.2, 2);
\draw[dashed, very thick] (-.5, -.5) -- (-.5, 2.5);
\draw[dashed, very thick] (2, -.5) -- (2, 2.5);
\draw[dotted, ultra thick, red] (1, 3) -- (3, 1);
\node[anchor=south east, color=red] at (1.1,2.9) (coord) {$\cC$};
\draw[dotted, ultra thick, red] (1.5,-0.5 ) -- (-0.5, 1.5);
\fill[pattern=north east lines, pattern color=gray] (.5, .5) -- (-.5,1.5) -- (-.5,2) -- (2,2) -- (2, .5) -- (.5,.5) -- cycle;%
\end{tikzpicture}
}
\caption{Results for Numerical Case I with {the} \priorB. The subfigures illustrate different scenarios for the prior belief: (a) it is a valid feasible set, (b) it is not a valid feasible set but a subset of known constraints  {$\cS$}, and (c) it is neither a valid feasible set nor a subset of {$\cS$}.   
}
\label{fig:priorB}
\end{figure}

\begin{figure}[htbp] 
\centering
\subfigure[\zeroObj \label{fig:noPriorB:a}]{   
\begin{tikzpicture} [thick, scale=1]
\draw[->, >=stealth] (-1,0) -- (3,0);
\draw[->, >=stealth] (0,-.5) -- (0,3);
\draw [fill, color=black] (1.5,1.5) circle [radius=0.07];
\draw [fill, color=black] (1,1) circle [radius=0.07];
\draw [fill, color=black] (2,1) circle [radius=0.07];
\draw [fill, color=black] (1,2) circle [radius=0.07];
\draw [fill, color=red] (2,2) circle [radius=0.07];
\node[anchor=south, color=red] at (2.2 , 2) (coord) {$\bx^0$};
\draw[<-] (3.2,2.2) coordinate -- (3.5,2.5) node[anchor= north west] {$\bc$};
\draw[-] (3.6, 2.4) -- (3.4, 2.6);
\draw[dashed,  very thick] (1, 3) -- (3, 1);;
\node[anchor=west, color=black] at (3 , 1) (coord) {\scriptsize ($4\times$)};
\draw[dotted,  ultra thick, red] (1, 3) -- (3, 1);
\node[anchor=south east, color=red] at (1.1,2.9) (coord) {$\cC$};
\draw[dotted,  ultra thick, red] (1.5,-0.5 ) -- (-0.5, 1.5);
\fill[pattern=north east lines, pattern color=gray] (-0.5, 1.5)--(1,3) -- (3,1) -- (1.5,-0.5) -- cycle;
\end{tikzpicture}
}
\subfigure[\minDist \label{fig:noPriorB:b}]{   
\begin{tikzpicture} [thick, scale=1]
\draw[->, >=stealth] (-1,0) -- (3,0);
\draw[->, >=stealth] (0,-.5) -- (0,3);
\draw [fill, color=black] (1.5,1.5) circle [radius=0.07];
\draw [fill, color=black] (1,1) circle [radius=0.07];
\draw [fill, color=black] (2,1) circle [radius=0.07];
\draw [fill, color=black] (1,2) circle [radius=0.07];
\draw [fill, color=red] (2,2) circle [radius=0.07];
\node[anchor=south, color=red] at (2.2 , 2) (coord) {$\bx^0$};
\draw[<-] (3.2,2.2) coordinate -- (3.5,2.5) node[anchor= north west] {$\bc$};
\draw[-] (3.6, 2.4) -- (3.4, 2.6);
\draw[dashed, very thick] (-0.5, 1) -- (3.5, 1); 
\node[anchor=west, color=black] at (3.5 , 1) (coord) {\scriptsize ($4\times$)};
\draw[dotted,  ultra thick, red] (1, 3) -- (3, 1);
\node[anchor=south east, color=red] at (1.1,2.9) (coord) {$\cC$};
\draw[dotted,  ultra thick, red] (1.5,-0.5 ) -- (-0.5, 1.5);
\fill[pattern=north east lines, pattern color=gray] (-0.5, 1.5)--(1,3) -- (3,1) -- (0,1) -- cycle;
\end{tikzpicture}
} 

\subfigure[\fairness \label{fig:noPriorB:c}]{   
\begin{tikzpicture} [thick, scale=1]
\draw[->, >=stealth] (-1,0) -- (3,0);
\draw[->, >=stealth] (0,-.5) -- (0,3);
\draw [fill, color=black] (1.5,1.5) circle [radius=0.07];
\draw [fill, color=black] (1,1) circle [radius=0.07];
\draw [fill, color=black] (2,1) circle [radius=0.07];
\draw [fill, color=black] (1,2) circle [radius=0.07];
\draw [fill, color=red] (2,2) circle [radius=0.07];
\node[anchor=south, color=red] at (2.2 , 2) (coord) {$\bx^0$};
\draw[<-] (3.2,2.2) coordinate -- (3.5,2.5) node[anchor= north west] {$\bc$};
\draw[-] (3.6, 2.4) -- (3.4, 2.6);
\draw[dashed, very thick] (.5, 3.5) -- (3, 1);
\draw[dashed, very thick] (1, 3.5) -- (1,.5);
\draw[dashed, very thick] (2, -.5) -- (2, 2.5);
\draw[dashed, very thick] (.5, 1.5) -- (2.5, -.5);
\draw[dotted, ultra thick, red] (1, 3) -- (3, 1);
\node[anchor=south west, color=red] at (1,2.9) (coord) {$\cC$};
\draw[dotted,  ultra thick, red] (1.5,-0.5 ) -- (-0.5, 1.5);
\fill[pattern=north east lines, pattern color=gray] (.5, .5) -- (1,3) -- (2, 2) -- (2, 0) -- (1,1) -- (1,3);
\end{tikzpicture}
}
\subfigure[\minmin \label{fig:noPriorB:d}]{   
\begin{tikzpicture} [thick, scale=1]
\draw[->, >=stealth] (-1,0) -- (3,0);
\draw[->, >=stealth] (0,-.5) -- (0,3);
\draw [fill, color=black] (1.5,1.5) circle [radius=0.07];
\draw [fill, color=black] (1,1) circle [radius=0.07];
\draw [fill, color=black] (2,1) circle [radius=0.07];
\draw [fill, color=black] (1,2) circle [radius=0.07];
\draw [fill, color=red] (2,2) circle [radius=0.07];
\node[anchor=south, color=red] at (2.2 , 2) (coord) {$\bx^0$};
\draw[<-] (3.2,2.2) coordinate -- (3.5,2.5) node[anchor= north west] {$\bc$};
\draw[-] (3.6, 2.4) -- (3.4, 2.6);
\draw[dashed, very thick] (-0.5, 1) -- (3.5, 1);
\draw[dashed, very thick] (-0.5, 2) -- (3.5, 2);
\draw[dashed, very thick] (1,-0.5 ) -- (1, 2.5);
\node[anchor=west, color=black] at (3.5 , 1) (coord) {\scriptsize ($2\times$)};
\draw[dotted, ultra thick, red] (1.5, 2.5) -- (3.5, .5);
\node[anchor=south east, color=red] at (1.6, 2.4) (coord) {$\cC$};
\draw[dotted, ultra thick, red] (1.5,-0.5 ) -- (-0.5, 1.5);
\fill[pattern=north east lines, pattern color=gray] 
(1, 1)--(1,2) -- (2,2) -- (3,1) --  cycle;
\end{tikzpicture}
} 
\caption{Results for Numerical Case I with different loss functions.}
\label{fig:noPriorB}
\end{figure}
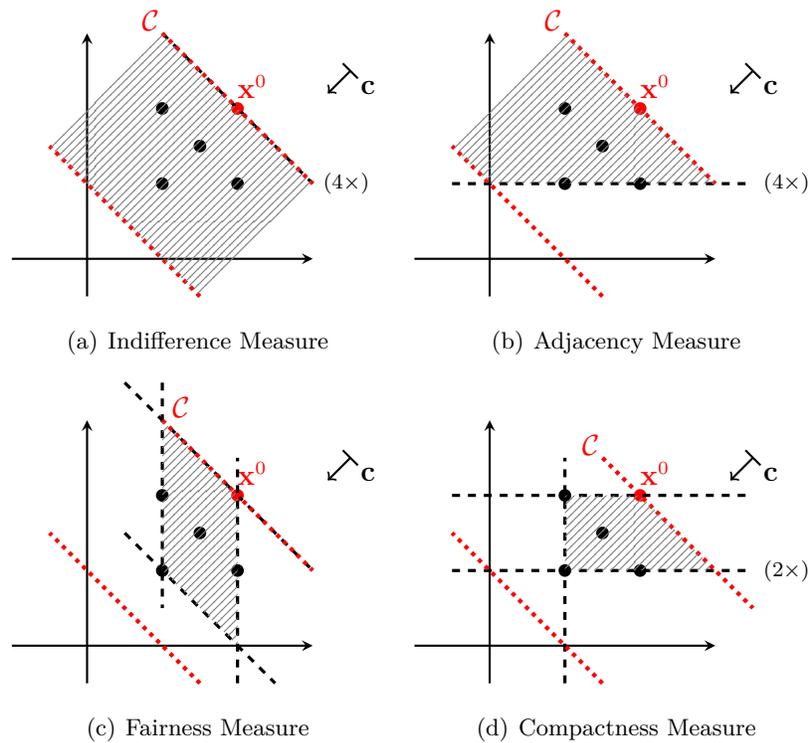

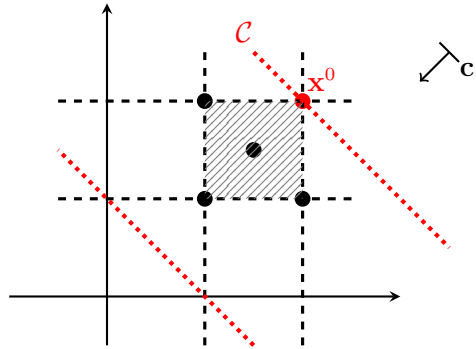
\begin{figure}[htbp]
    \centering
\begin{tikzpicture} [thick, scale=1.3]
\draw[->, >=stealth] (-1,0) -- (3,0);
\draw[->, >=stealth] (0,-.5) -- (0,3);
\draw [fill, color=black] (1.5,1.5) circle [radius=0.07];
\draw [fill, color=black] (1,1) circle [radius=0.07];
\draw [fill, color=black] (2,1) circle [radius=0.07];
\draw [fill, color=black] (1,2) circle [radius=0.07];
\draw [fill, color=red] (2,2) circle [radius=0.07];
\node[anchor=south, color=red] at (2.2 , 2) (coord) {$\bx^0$};
\draw[<-] (3.2,2.2) coordinate -- (3.5,2.5) node[anchor= north west] {$\bc$};
\draw[-] (3.6, 2.4) -- (3.4, 2.6);
\draw[dashed, very thick] (-.5, 1) -- (2.5, 1);
\draw[dashed, very thick] (1, -.5) -- (1,2.5);
\draw[dashed, very thick] (-.5, 2) -- (2.5, 2);
\draw[dashed, very thick] (2, -.5) -- (2, 2.5);
\draw[dotted, ultra thick, red] (1.5, 2.5) -- (3.5, .5);
\node[anchor=south east, color=red] at (1.6,2.5) (coord) {$\cC$};
\draw[dotted, ultra thick, red] (1.5,-0.5 ) -- (-0.5, 1.5);
\fill[pattern=north east lines, pattern color=gray] (.5, .5) -- (1,1) -- (1, 2) -- (2, 2) -- (2,1) -- (1,1);
\end{tikzpicture}
\caption{Result for Numerical Case I with a combined loss function, where the \fairness~is the primary objective and the \minDist~is the secondary objective. }
\label{fig:secondObj}
\end{figure}

For this example, we present the results for a loss function with \priorB, the four loss function when no prior belief is provided, and an example of the combined loss functions, in Figures~\ref{fig:priorB},~\ref{fig:noPriorB}, and~\ref{fig:secondObj}, respectively. In these figures, black dots denote the given observations $\bx^k$, $\forall k \in \cK$, and the preferred observation ($\bx^0$) is highlighted in red. The blue solid lines are the hyperplanes corresponding to the given prior belief parameters ($\barhatA$), the dotted red lines represent the known constraints ($\cS$), and the dashed black lines demonstrate the constraints found by the inverse optimization model. The resulting imputed feasible set of $\MIO$ (including the known constraints) is marked as a shaded area.

Figure~\ref{fig:priorB} shows the results under the \priorB~for three different possible feasible regions as the prior belief ($\barhatA$). With the \priorB, the goal of the inverse optimization model is to minimally perturb  $\barhatA$ to ensure all the observations are feasible and $\bx^0$ is optimal. In Figure~\ref{fig:priorB:a}, the given prior belief $\barhatA$ is a valid feasible set, and as shown in Proposition~\ref{prop:prior2}, the optimal solution $\barA$ is the same as the prior belief $\barhatA$. On the contrary, in Figures~\ref{fig:priorB:b} and~\ref{fig:priorB:c}, the given prior beliefs are not valid feasible sets. In Figure~\ref{fig:priorB:b}, $\barhatA \subseteq \cS$, while in Figure~\ref{fig:priorB:c},  a part of the prior belief is infeasible for the known constraints and hence, $\barhatA \not\subseteq \cS$. In both cases, the solution $\barA$ is a minimally perturbed $\barhatA$ that makes all the observations feasible. The resulting imputed feasible set (\ie, the shaded area) is derived from $\barA$ and is a subset of $\cS$. 

Figure~\ref{fig:priorB} confirms that the \priorB~heavily relies on the quality of the prior belief on the constraint parameters. Particularly, when the prior belief is an unreasonably large valid feasible set far from the observations (\eg, Figure~\ref{fig:priorB:a}, or two of the constraints in Figure~\ref{fig:priorB:c}) the inverse problem will always return the prior belief as the optimal solution. While this measure is the most commonly-used objective function of inverse problems in the literature, the results obtained may not be reliable if a high-quality prior belief does not exist. 

Figure~\ref{fig:noPriorB} illustrates the results for the other loss functions as defined in Section~\ref{sec:measures}. Figure~\ref{fig:noPriorB:a} shows the results for the \zeroObj~which has infinitely many optimal solutions. In our results, the four inferred constraints happened to be the same and equal to the half-space $\cC$. Hence, the resulting inferred feasible region is equivalent to the set of known constraints~$\cS$, which is the largest possible imputed feasible set (Remark~\ref{rem:cs}). Figures~\ref{fig:noPriorB:b}, \ref{fig:noPriorB:c}, and \ref{fig:noPriorB:d} show the results for the \minDist, the \fairness, and the \minmin, respectively. In Figure~\ref{fig:noPriorB:b}, the inferred constraints are four identical lines that pass through observations $(1,1)$ and $(2,1)$. On the contrary, Figure~\ref{fig:noPriorB:c}, shows that employing the \fairness~in the objective function results in finding four distinct constraints. For this example, the constraints have the same total distance from all observations and are hence, distributed fairly. As speculated, this measure provides a confined feasible set for the $\FO$ problem. Finally, Figure~\ref{fig:noPriorB:d} illustrates the results for the \minmin which ensures that each observation is close to some constraint. In this case, each constraint passes through two of the observations, and due to symmetry, two of the constraints are identical.

Lastly, we used a combined loss function to provide additional control over the properties of the imputed feasible set as shown in Figure \ref{fig:secondObj}. In this example, we first imposed the \fairness~to encourage similar total distances across different constraints, and then applied the \minDist~as a secondary objective. In other words, we search among those solutions with the optimal \fairness~value that also have the minimum total distance between constraints and observations. As Figure~\ref{fig:secondObj} illustrates, the imputed feasible set is the same as the convex hull $\cH$ in this case, which is the smallest possible imputed feasible set for $\FO$ (Lemma~\ref{lem:subset}).

\subsection{Numerical Case II}\label{sec:NumCase2}
In this section, we test our approach on a relatively larger numerical example. This example considers 19 observations with $\bx^0 = (1,1)$, two known constraints (the first one being the half-space $\cC$), and 6 unknown constraints. The details of this numerical example are summarized in Table~\ref{tab:case2}. Given the larger size of this example, the nonlinear solver {MINOS} was not able to find the global optimal solutions for most instances of the problem. Hence, we were only able to solve this example using the $\eMIO$ formulation, which further illustrates the importance and advantage of the proposed $\eMIO$ formulation in solving larger instances to optimality. This advantage is more pronounced when the normalization constraint~\eqref{eq:MIPnorm} is linear and the $\eMIO$ problem becomes a linearly-constrained optimization model. 

\begin{table}[htbp]
    \centering
 \begin{tabular}{p{.28\linewidth} p{.6\linewidth}}
    \toprule
    {\bf Description} & {\bf Value(s)} \\
    \midrule 
    Cost vector $\bc$ & $(1, 1)$ \\
    \midrule
    Observations $\bx^0; \bx^k \qquad$  & $\bf{(1,1)}$;  $(2,1)$, $(4,2)$, $(4,5)$, $(3,6)$, $(2,4)$, $(3,4)$, $(3,2)$, $(4,3)$, $(1,3)$, $(2,2.5)$, $(1,5)$, $(5,2.5)$, $(5,4)$, $(2.7,3.2)$, $(2.3,4.7)$, $(1.4,4.8)$, $(3.8,4.3)$, $(4.8,3.3)$ \\
    \midrule
    Known constraints   
                        & $0.5 x_1 + 0.5 x_2 \geq 1 \quad$  (half-space $\cC$) \\ 
                        &  $-x_1  \geq -5 $ \\
    \midrule
    Unknown constraints & 6 constraints\\
    \bottomrule
\end{tabular}
    \caption{Numerical Case II}
    \label{tab:case2}
\end{table}

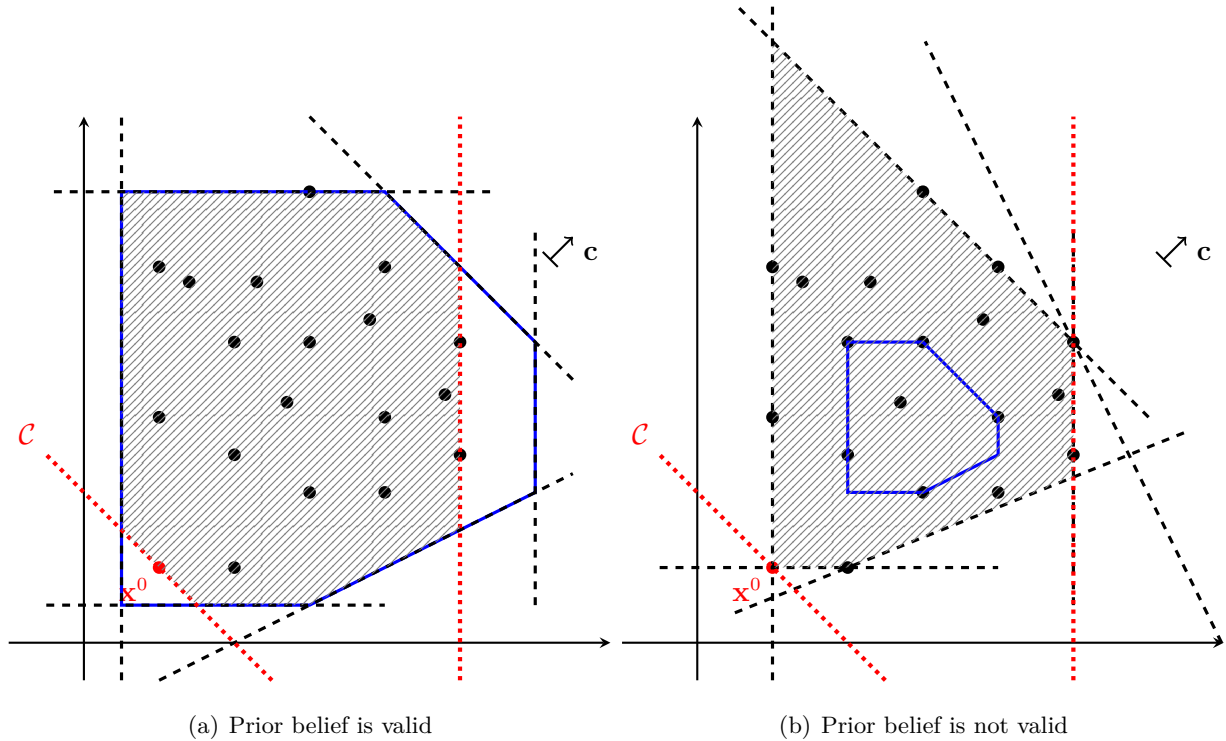
\begin{figure}[htbp] \vspace{-2em}
\centering
\subfigure[Prior belief is valid \label{fig:largefigPrior:a}]
{%
\begin{tikzpicture}[thick, scale=1] 
\draw[->, >=stealth] (-1,0) -- (7,0);
\draw[->, >=stealth] (0,-.5) -- (0,7);
\draw [fill, color=black] (2,1) circle [radius=0.07];
\draw [fill, color=black] (4,2) circle [radius=0.07];
\draw [fill, color=black] (5,2.5) circle [radius=0.07];
\draw [fill, color=black] (5,4) circle [radius=0.07];
\draw [fill, color=black] (4,5) circle [radius=0.07];
\draw [fill, color=black] (3,6) circle [radius=0.07];
\draw [fill, color=black] (1,5) circle [radius=0.07];
\draw [fill, color=black] (2,4) circle [radius=0.07];
\draw [fill, color=black] (3,4) circle [radius=0.07];
\draw [fill, color=black] (3,2) circle [radius=0.07];
\draw [fill, color=black] (4,3) circle [radius=0.07];
\draw [fill, color=black] (1,3) circle [radius=0.07];
\draw [fill, color=black] (2,2.5) circle [radius=0.07];
\draw [fill, color=black] (2.7,3.2) circle [radius=0.07];
\draw [fill, color=black] (2.3,4.8) circle [radius=0.07];
\draw [fill, color=black] (1.4,4.8) circle [radius=0.07];
\draw [fill, color=black] (3.8,4.3) circle [radius=0.07];
\draw [fill, color=black] (4.8,3.3) circle [radius=0.07];
%
\draw [fill, color=red] (1,1) circle [radius=0.07]; 
\node[anchor=north east, color=red] at (1 , 1) (coord) {$\bx^0$};
%
\draw[->] (6.2,5.1) coordinate -- (6.5,5.4) node[anchor=north west] {$\bc$};
\draw[-] (6.1,5.20) -- (6.30, 5.0);
\draw[blue, very thick] (0.5,0.5) -- (0.5,6) -- (4,6) -- (6,4) -- (6,2) -- (3, 0.5) -- cycle ;
\draw[dashed, very thick] (0.5, -0.5) -- (0.5, 7);
\draw[dashed, very thick] (-0.5, 0.5) -- (4,0.5);
\draw[dashed, very thick] (1,-0.5) -- (6.5,2.25);
\draw[dashed, very thick] (6,0.5) -- (6,5.5);
\draw[dashed, very thick] (6.5,3.5) -- (3,7);
\draw[dashed, very thick] (5.4,6) -- (-0.5,6);
\draw[dotted, red, ultra thick] (2.5,-0.5) -- (-0.5,2.5);
\node[anchor=south east, color=red] at (-0.5, 2.5) (coord) {$\cC$};
\draw[dotted, red, ultra thick] (5,-0.5) -- (5,7);
\fill[pattern=north east lines, pattern color=gray] (1.5,0.5) -- (0.5,1.5) --(0.5,6) -- (4,6) -- (5,5) -- (5,1.5) -- (3,0.5) -- cycle;
\end{tikzpicture}
}
\subfigure[Prior belief is not valid \label{fig:largefigPrior:b}]
{%
\begin{tikzpicture}[thick, scale=1] 
\draw[->, >=stealth] (-1,0) -- (7,0);
\draw[->, >=stealth] (0,-.5) -- (0,7);
\draw [fill, color=black] (2,1) circle [radius=0.07];
\draw [fill, color=black] (4,2) circle [radius=0.07];
\draw [fill, color=black] (5,2.5) circle [radius=0.07];
\draw [fill, color=black] (5,4) circle [radius=0.07];
\draw [fill, color=black] (4,5) circle [radius=0.07];
\draw [fill, color=black] (3,6) circle [radius=0.07];
\draw [fill, color=black] (1,5) circle [radius=0.07];
\draw [fill, color=black] (2,4) circle [radius=0.07];
\draw [fill, color=black] (3,4) circle [radius=0.07];
\draw [fill, color=black] (3,2) circle [radius=0.07];
\draw [fill, color=black] (4,3) circle [radius=0.07];
\draw [fill, color=black] (1,3) circle [radius=0.07];
\draw [fill, color=black] (2,2.5) circle [radius=0.07];
\draw [fill, color=black] (2.7,3.2) circle [radius=0.07];
\draw [fill, color=black] (2.3,4.8) circle [radius=0.07];
\draw [fill, color=black] (1.4,4.8) circle [radius=0.07];
\draw [fill, color=black] (3.8,4.3) circle [radius=0.07];
\draw [fill, color=black] (4.8,3.3) circle [radius=0.07];
%
\draw [fill, color=red] (1,1) circle [radius=0.07]; 
\node[anchor=north east, color=red] at (1 , 1) (coord) {$\bx^0$};
%
\draw[->] (6.2,5.1) coordinate -- (6.5,5.4) node[anchor=north west] {$\bc$};
\draw[-] (6.1,5.20) -- (6.30, 5.0);
\draw[blue, very thick] (2,2) -- (2,4) -- (3,4) -- (4,3) -- (4,2.5) -- (3,2) -- cycle;
\draw[dashed, very thick] (1, -0.5) -- (1, 8.5);
\draw[dashed, very thick] (-0.5, 1) -- (4,1);
\draw[dashed, very thick] (.5,0.398) -- (6.5,2.805);
\draw[dashed, very thick] (5,0.5) -- (5,5.5);
\draw[dashed, very thick] (6,3) -- (.5,8.5);
 \draw[dashed, very thick] (6.97, 0) -- (3.03,8);
\draw[dotted, red, ultra thick] (2.5,-0.5) -- (-0.5,2.5);
\node[anchor=south east, color=red] at (-0.5, 2.5) (coord) {$\cC$};
\draw[dotted, red, ultra thick] (5,-0.5) -- (5,7);
\fill[pattern=north east lines, pattern color=gray] (1,1) -- (1,8) -- (5,4) -- (5,2.204) --(2,1) -- cycle;
\end{tikzpicture}
}
\caption{Results for Numerical Case II with the \priorB. (a) The prior belief $\barhatA$ is a valid feasible set but $\barhatA\not\subseteq\cS$. (b) The prior belief is not a valid feasible set but $\barhatA\subseteq \cS$. 
}
\label{fig:largefigPrior}
\end{figure}

\begin{figure}[htbp] \vspace{-1.5em}
\centering
\subfigure[\zeroObj \label{fig:largefigNoPrior:a}]
{%
\begin{tikzpicture}[thick, scale=0.9] 
\draw[->, >=stealth] (-1,0) -- (7,0);
\draw[->, >=stealth] (0,-.5) -- (0,7);
\draw [fill, color=black] (2,1) circle [radius=0.07];
\draw [fill, color=black] (4,2) circle [radius=0.07];
\draw [fill, color=black] (5,2.5) circle [radius=0.07];
\draw [fill, color=black] (5,4) circle [radius=0.07];
\draw [fill, color=black] (4,5) circle [radius=0.07];
\draw [fill, color=black] (3,6) circle [radius=0.07];
\draw [fill, color=black] (1,5) circle [radius=0.07];
\draw [fill, color=black] (2,4) circle [radius=0.07];
\draw [fill, color=black] (3,4) circle [radius=0.07];
\draw [fill, color=black] (3,2) circle [radius=0.07];
\draw [fill, color=black] (4,3) circle [radius=0.07];
\draw [fill, color=black] (1,3) circle [radius=0.07];
\draw [fill, color=black] (2,2.5) circle [radius=0.07];
\draw [fill, color=black] (2.7,3.2) circle [radius=0.07];
\draw [fill, color=black] (2.3,4.8) circle [radius=0.07];
\draw [fill, color=black] (1.4,4.8) circle [radius=0.07];
\draw [fill, color=black] (3.8,4.3) circle [radius=0.07];
\draw [fill, color=black] (4.8,3.3) circle [radius=0.07];
%
\draw [fill, color=red] (1,1) circle [radius=0.07]; 
\node[anchor=north east, color=red] at (1 , 1) (coord) {$\bx^0$};
%
\draw[->] (6.2,5.1) coordinate -- (6.5,5.4) node[anchor=north west] {$\bc$};
\draw[-] (6.1,5.20) -- (6.30, 5.0);
\draw[dashed, very thick] (1,-0.5) -- (1,7);
\node[anchor=south, color=black] at (1 , 7) (coord) {\scriptsize $(6\times)$};
\draw[dotted, red, ultra thick] (2.5,-0.5) -- (-0.5,2.5);
\node[anchor=south east, color=red] at (-0.5, 2.5) (coord) {$\cC$};
\draw[dotted, red, ultra thick] (5,-0.5) -- (5,7);
\fill[pattern=north east lines, pattern color=gray] (2,0) --(1,1) -- (1,7) -- (5,7) -- (5,-0.5) -- (2.5,-0.5) -- cycle;
\end{tikzpicture}
}
\subfigure[\minDist \label{fig:largefigNoPrior:b}]
{%
\begin{tikzpicture}[thick, scale=0.9] 
\draw[->, >=stealth] (-1,0) -- (7,0);
\draw[->, >=stealth] (0,-.5) -- (0,7);
\draw [fill, color=black] (2,1) circle [radius=0.07];
\draw [fill, color=black] (4,2) circle [radius=0.07];
\draw [fill, color=black] (5,2.5) circle [radius=0.07];
\draw [fill, color=black] (5,4) circle [radius=0.07];
\draw [fill, color=black] (4,5) circle [radius=0.07];
\draw [fill, color=black] (3,6) circle [radius=0.07];
\draw [fill, color=black] (1,5) circle [radius=0.07];
\draw [fill, color=black] (2,4) circle [radius=0.07];
\draw [fill, color=black] (3,4) circle [radius=0.07];
\draw [fill, color=black] (3,2) circle [radius=0.07];
\draw [fill, color=black] (4,3) circle [radius=0.07];
\draw [fill, color=black] (1,3) circle [radius=0.07];
\draw [fill, color=black] (2,2.5) circle [radius=0.07];
\draw [fill, color=black] (2.7,3.2) circle [radius=0.07];
\draw [fill, color=black] (2.3,4.8) circle [radius=0.07];
\draw [fill, color=black] (1.4,4.8) circle [radius=0.07];
\draw [fill, color=black] (3.8,4.3) circle [radius=0.07];
\draw [fill, color=black] (4.8,3.3) circle [radius=0.07];
%
\draw [fill, color=red] (1,1) circle [radius=0.07]; 
\node[anchor=north east, color=red] at (1 , 1) (coord) {$\bx^0$};
%
\draw[->] (6.2,5.1) coordinate -- (6.5,5.4) node[anchor=north west] {$\bc$};
\draw[-] (6.1,5.20) -- (6.30, 5.0);
\draw[dashed, very thick] (2,7) -- (6,3);
\node[anchor=south, color=black] at (2.2 , 7) (coord) {\scriptsize $(6\times)$};
\draw[dotted, red, ultra thick] (2.5,-0.5) -- (-0.5,2.5);
\node[anchor=south east, color=red] at (-0.5, 2.5) (coord) {$\cC$};
\draw[dotted, red, ultra thick] (5,-0.5) -- (5,7);
\fill[pattern=north east lines, pattern color=gray] (2,0) -- (-.5,2.5) -- (-.5,7) --(2,7)-- (5,4) -- (5,-0.5) -- (2.5, -0.5) -- cycle;
\end{tikzpicture}
}

\subfigure[\fairness \label{fig:largefigNoPrior:c}]
{%
\begin{tikzpicture}[thick, scale=0.9] 
\draw[->, >=stealth] (-1,0) -- (7,0);
\draw[->, >=stealth] (0,-.5) -- (0,7);
\draw [fill, color=black] (2,1) circle [radius=0.07];
\draw [fill, color=black] (4,2) circle [radius=0.07];
\draw [fill, color=black] (5,2.5) circle [radius=0.07];
\draw [fill, color=black] (5,4) circle [radius=0.07];
\draw [fill, color=black] (4,5) circle [radius=0.07];
\draw [fill, color=black] (3,6) circle [radius=0.07];
\draw [fill, color=black] (1,5) circle [radius=0.07];
\draw [fill, color=black] (2,4) circle [radius=0.07];
\draw [fill, color=black] (3,4) circle [radius=0.07];
\draw [fill, color=black] (3,2) circle [radius=0.07];
\draw [fill, color=black] (4,3) circle [radius=0.07];
\draw [fill, color=black] (1,3) circle [radius=0.07];
\draw [fill, color=black] (2,2.5) circle [radius=0.07];
\draw [fill, color=black] (2.7,3.2) circle [radius=0.07];
\draw [fill, color=black] (2.3,4.8) circle [radius=0.07];
\draw [fill, color=black] (1.4,4.8) circle [radius=0.07];
\draw [fill, color=black] (3.8,4.3) circle [radius=0.07];
\draw [fill, color=black] (4.8,3.3) circle [radius=0.07];
%
\draw [fill, color=red] (1,1) circle [radius=0.07]; 
\node[anchor=north east, color=red] at (1 , 1) (coord) {$\bx^0$};
%
\draw[->] (6.2,5.1) coordinate -- (6.5,5.4) node[anchor=north west] {$\bc$};
\draw[-] (6.1,5.20) -- (6.30, 5.0);
\draw[dashed, very thick] (2.938,-0.938) -- (5.96,4.1);
\draw[dashed, very thick] (1,-0.5) -- (1,7) ;
\draw[dashed, very thick] (0.2,1)--(1.6,8);
\node[anchor=west, color=black] at (1.6,8) (coord) {\scriptsize $(2\times)$};
\draw[dashed, very thick] (1,8)--(6,3);
\node[anchor=north  , color=black] at (6,3) (coord) {\scriptsize $(2\times)$};
\draw[dotted, red, ultra thick] (2.5,-0.5) -- (-0.5,2.5);
\node[anchor=south east, color=red] at (-0.5, 2.5) (coord) {$\cC$};
\draw[dotted, red, ultra thick] (5,-0.5) -- (5,7);
\fill[pattern=north east lines, pattern color=gray] (2.938,-0.938) --(1,1)-- (1,5)--(1.5,7.5)--(5,4) -- (5,2.5) -- (3.5, 0)-- cycle;
\end{tikzpicture}
}
\subfigure[\minmin \label{fig:largefigNoPrior:d}]
{%
\begin{tikzpicture}[thick, scale=0.9] 
\draw[->, >=stealth] (-1,0) -- (7,0);
\draw[->, >=stealth] (0,-.5) -- (0,7);
\draw [fill, color=black] (2,1) circle [radius=0.07];
\draw [fill, color=black] (4,2) circle [radius=0.07];
\draw [fill, color=black] (5,2.5) circle [radius=0.07];
\draw [fill, color=black] (5,4) circle [radius=0.07];
\draw [fill, color=black] (4,5) circle [radius=0.07];
\draw [fill, color=black] (3,6) circle [radius=0.07];
\draw [fill, color=black] (1,5) circle [radius=0.07];
\draw [fill, color=black] (2,4) circle [radius=0.07];
\draw [fill, color=black] (3,4) circle [radius=0.07];
\draw [fill, color=black] (3,2) circle [radius=0.07];
\draw [fill, color=black] (4,3) circle [radius=0.07];
\draw [fill, color=black] (1,3) circle [radius=0.07];
\draw [fill, color=black] (2,2.5) circle [radius=0.07];
\draw [fill, color=black] (2.7,3.2) circle [radius=0.07];
\draw [fill, color=black] (2.3,4.8) circle [radius=0.07];
\draw [fill, color=black] (1.4,4.8) circle [radius=0.07];
\draw [fill, color=black] (3.8,4.3) circle [radius=0.07];
\draw [fill, color=black] (4.8,3.3) circle [radius=0.07];
%
\draw [fill, color=red] (1,1) circle [radius=0.07]; 
\node[anchor=north east, color=red] at (1 , 1) (coord) {$\bx^0$};
\draw[->] (6.2,5.1) coordinate -- (6.5,5.4) node[anchor=north west] {$\bc$};
\draw[-] (6.1,5.20) -- (6.30, 5.0);
\draw[dashed, very thick] (-1,-0.5) -- (6.5,3.25);
\node[anchor=south, color=black] at (6.5 , 3.25) (coord) {\scriptsize $(2\times)$};
\draw[dashed, very thick] (1,-0.5) -- (1,8.5);
\node[anchor=west, color=black] at (1,8.4) (coord) {\scriptsize $(2\times)$};
\draw[dashed, very thick] (0.5,8.5) -- (6.5,2.5);
\draw[dashed, very thick] (5,-0.5) -- (5,7);
\draw[dotted, red, ultra thick] (2.5,-0.5) -- (-0.5,2.5);
\node[anchor=south east, color=red] at (-0.5, 2.5) (coord) {$\cC$};
\draw[dotted, red, ultra thick] (5,-0.5) -- (5,7);
\fill[pattern=north east lines, pattern color=gray] (1.333, 0.667) -- (1,1) -- (1,8) -- (5,4) -- (5,2.5) -- cycle;
\end{tikzpicture}
}
\caption{Results for Numerical Case II with different loss functions.}
\label{fig:largefigNoPrior}
\end{figure}
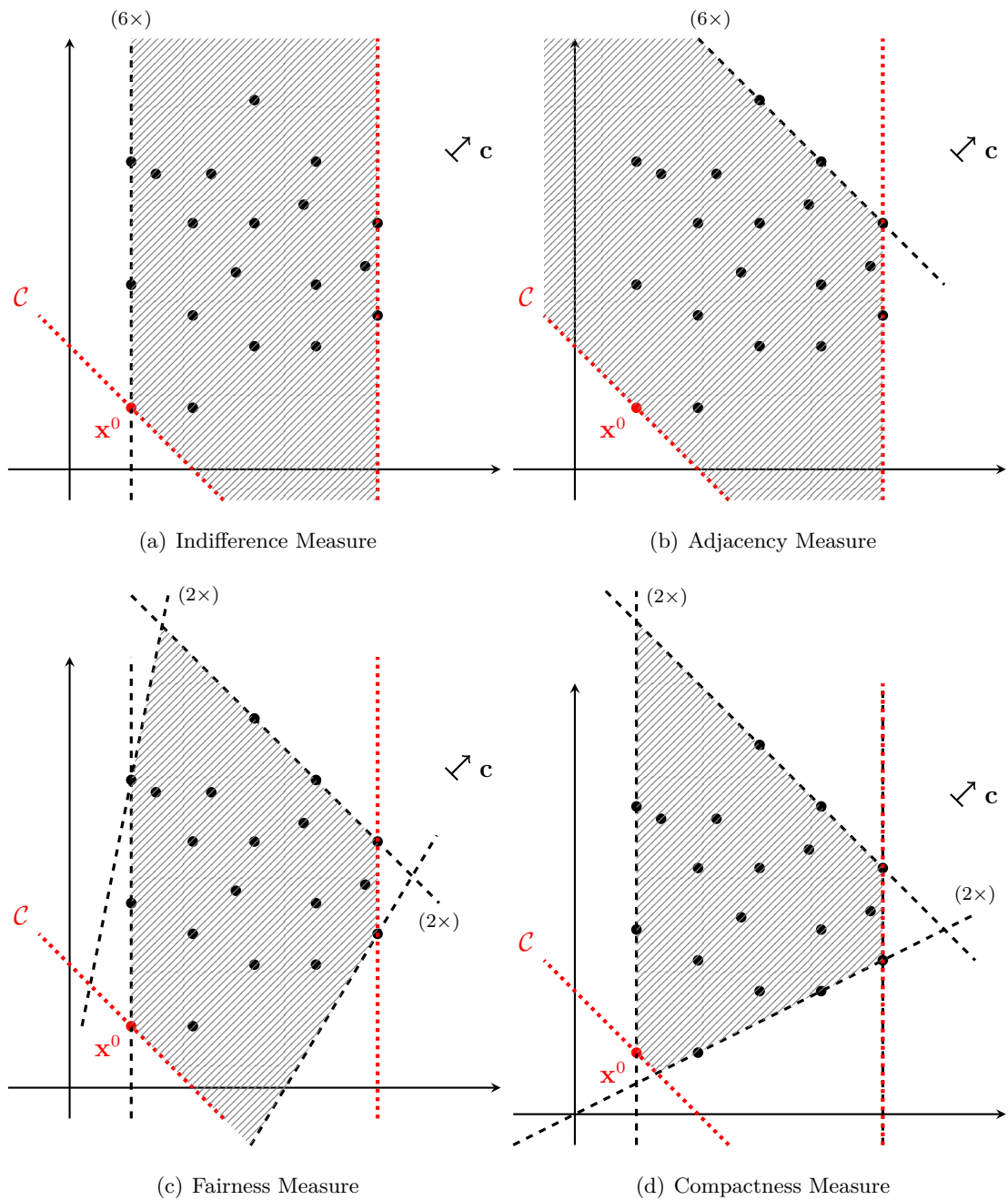

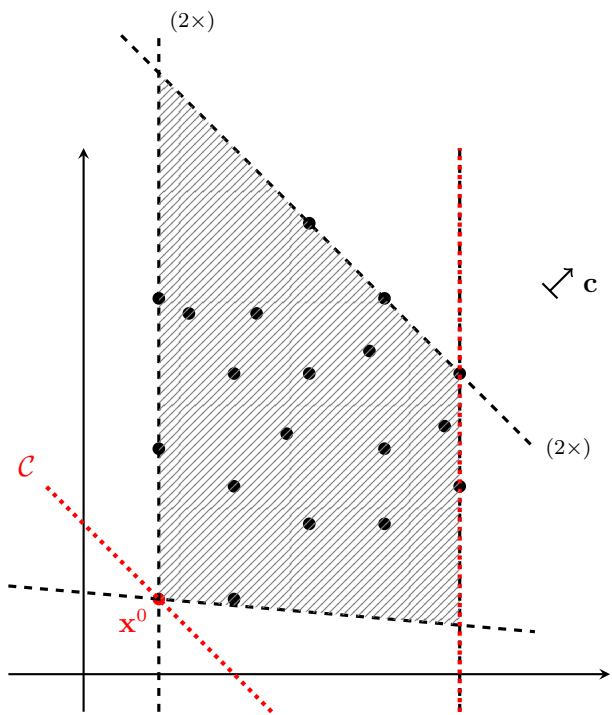
\begin{figure}[htbp] 
    \centering
\begin{tikzpicture}[thick, scale=1] 
\draw[->, >=stealth] (-1,0) -- (7,0);
\draw[->, >=stealth] (0,-.5) -- (0,7);
\draw [fill, color=black] (2,1) circle [radius=0.07];
\draw [fill, color=black] (4,2) circle [radius=0.07];
\draw [fill, color=black] (5,2.5) circle [radius=0.07];
\draw [fill, color=black] (5,4) circle [radius=0.07];
\draw [fill, color=black] (4,5) circle [radius=0.07];
\draw [fill, color=black] (3,6) circle [radius=0.07];
\draw [fill, color=black] (1,5) circle [radius=0.07];
\draw [fill, color=black] (2,4) circle [radius=0.07];
\draw [fill, color=black] (3,4) circle [radius=0.07];
\draw [fill, color=black] (3,2) circle [radius=0.07];
\draw [fill, color=black] (4,3) circle [radius=0.07];
\draw [fill, color=black] (1,3) circle [radius=0.07];
\draw [fill, color=black] (2,2.5) circle [radius=0.07];
\draw [fill, color=black] (2.7,3.2) circle [radius=0.07];
\draw [fill, color=black] (2.3,4.8) circle [radius=0.07];
\draw [fill, color=black] (1.4,4.8) circle [radius=0.07];
\draw [fill, color=black] (3.8,4.3) circle [radius=0.07];
\draw [fill, color=black] (4.8,3.3) circle [radius=0.07];
%
\draw [fill, color=red] (1,1) circle [radius=0.07]; 
\node[anchor=north east, color=red] at (1 , 1) (coord) {$\bx^0$};
%
\draw[->] (6.2,5.1) coordinate -- (6.5,5.4) node[anchor=north west] {$\bc$};
\draw[-] (6.1,5.20) -- (6.30, 5.0);
\draw[dashed, very thick] (1,-0.5) -- (1,8.5);
\node[anchor=south  west, color=black] at (1,8.4) (coord) {\scriptsize $(2\times)$};
\draw[dashed, very thick] (0.5,8.5) -- (6,3);
 \node[anchor=west, color=black] at (6 , 3) (coord) {\scriptsize $(2\times)$};
\draw[dashed, very thick] (-1,1.174) -- (6,0.565);
\draw[dashed, very thick] (5,-0.5) -- (5,7);
\draw[dotted, red, ultra thick] (2.5,-0.5) -- (-0.5,2.5);
\node[anchor=south east, color=red] at (-0.5, 2.5) (coord) {$\cC$};
\draw[dotted, red, ultra thick] (5,-0.5) -- (5,7);
 \fill[pattern=north east lines, pattern color=gray] (1,1) -- (1,8) -- (5,4) -- (5,0.652) -- cycle;
\end{tikzpicture}
\caption{Result for Numerical Case II with a combined loss function, where the  \fairness~is the primary objective and the \minDist~is the secondary objective. }
    \label{fig:largefig-secondObj}
\end{figure}

Figure~\ref{fig:largefigPrior} illustrates the results for the \priorB. In Figure~\ref{fig:largefigPrior:a}, the prior belief is a valid feasible set, and the optimal solution $\barA$ is the same as the given prior belief $\barhatA$, as demonstrated by Proposition~\ref{prop:prior2}. Note that in this example, although $\barhatA$ is a valid feasible set, it does not satisfy the known constraints (\ie, $\barhatA \not\subseteq \cS$). Conversely, Figure~\ref{fig:largefigPrior:b} illustrates the case that $\barhatA$ is a subset of the known constraints but is not a valid feasible set. In this case, the prior belief is minimally expanded in order to include all observations. These results re-emphasize that the shape of the imputed feasible set is heavily affected by the quality of the prior belief on the constraint parameters. 

We next illustrate the results for the Numerical Case II with the remaining four loss functions that do not require a prior belief in Figure~\ref{fig:largefigNoPrior}. Analogous to Numerical Case I, Figure~\ref{fig:largefigNoPrior:a} shows the results for \zeroObj~which is a simple feasibility problem and results in the optimality of any feasible solution to $\MIO$. In our results, all six constraints equate a hyperplane that passes through three of the observations. Figure~\ref{fig:largefigNoPrior:b} shows the results for the \minDist~which, as expected, induces an unbounded imputed feasible set. In this case, the problem finds a constraint that has the minimum distance from all observations and chooses to repeat that constraint 6 times. 

Figures~\ref{fig:largefigNoPrior:c} and~\ref{fig:largefigNoPrior:d} demonstrate the results for the \fairness~and \minmin~as the loss function of $\MIO$, respectively. In both cases, the imputed feasible sets are bounded and $\bx^0$ is an extreme point. The \fairness~tries to spread the constraints around the observations and obtains a bounded imputed feasible set for this example. Finally, the \minmin~ensures that there is at least one constraint in the vicinity of each of the observations on the boundary of their convex hull results in an imputed feasible region that encapsulates all observations for this example.

As described in Section~\ref{sec:measures}, inverse models often have multiple solutions. Using combined loss functions, we can explore the multiple optimal solution space under one loss function and further tailor the desired characteristics of the imputed feasible set based on an alternate loss function. In Figure~\ref{fig:largefig-secondObj}, we use a combined loss function to first solve the inverse optimization problem with the \fairness~and then search among its multiple optimal solutions to find a set of constraints that minimizes the \minDist. As a result, we derive an imputed feasible set that not only scatters the constraints around the observations fairly but also minimizes the total distance of the constraints to all observations. In this example, the resulting imputed feasible region is visually tighter than that of the \fairness~alone.

{\color{myGreen}

\subsection{Numerical Case III} 
In this section, we validate our proposed methodology on a diet recommendation problem using a dataset of 100 observations of daily food intake choices. As briefly discussed in Section~\ref{sec:motivation}, the goal is to employ our $\MIO$ framework to impute the implicit constraints of a dieter based on the observations of past food choices. The implicit constraints are typically difficult to capture in regular forward settings. However, using our $\MIO$ approach, these additional constraints will help to identify diets that are more palatable to the user.

\begin{table}[htbp]
    \centering
 \begin{tabular}{p{.28\linewidth} p{.6\linewidth}}
    \toprule
    {\bf Description} & {\bf Value(s)} \\
    \midrule 
    Cost vector $\bc$ & ({\it a}\,) Maximize total protein intake \\
                      & ({\it b}\,) Minimize total sodium intake\\
    \midrule
    Observations &  100 Daily food intakes of 26 different food items \\
    \midrule
    Known constraints  
                        & 8 known constraints \& half-space $\cC$\\
    \it \footnotesize \hspace{0.9cm} Lower bounds: & \it \footnotesize Carbohydrates, Fiber, Calories \\
    \it \footnotesize \hspace{0.9cm} Upper bounds: & \it \footnotesize Fat, Sugar, Cholesterol, Calories, \# of servings \\ 
     \midrule
    Unknown constraints & 30 constraints\\
    \bottomrule
\end{tabular}
    \caption{Numerical Case III: Two objective functions were considered for a set of 100 observations on 26 food items. The set of constraints includes 8 known constraints, the half-space, and 30 unknown constraints.}
    \label{tab:case3}
\end{table}

In this case study, we consider a set of 100 observations of the daily food intake that are reported as the number of servings consumed per food per day of observation. This data is gathered from \citet{CSSEDietData} which is based on the National Health and Nutrition Examination Survey (NHANES) dietary data. Our data includes 26 food items, each of which were consumed at least once among the 100 daily observations. We consider a set of 8 known constraints on different nutrition values (\eg, lower bound on fiber, upper bound on cholesterol) that must always be met. To derive the known constraints, we consulted the guidelines provided by the~\citet{DietData_Const}. 

We tested the proposed $\MIO$ model with two different known objective functions: ({\it a}\,) maximizing daily protein intake and ({\it b}\,) minimizing daily sodium intake. For each objective, we found the preferred observation $\bx^0$ (\ie, the observation with the best objective value) and recovered a feasible region that made all 100 observations feasible and $\bx^0$ optimal for the corresponding forward problem. Table~\ref{tab:case3} shows a summary of the known constraints and the inverse problems, and Table~\ref{tab:case3_fooddetails} in the Appendix provides additional information about the observations. 
For each of the two objective functions, we solved the $\eMIO$ formulation to impute 30 additional constraints using a combined loss function with \fairness\ and \minmin\ as the primary and secondary objective functions, respectively. We then found optimal diets by solving the $\FO$ model twice: first with only the eight given constraints and then using both the known constraints and the imputed constraints by $\MIO$. 

\begin{figure}
\centering
    \includegraphics[width=.9\textwidth]{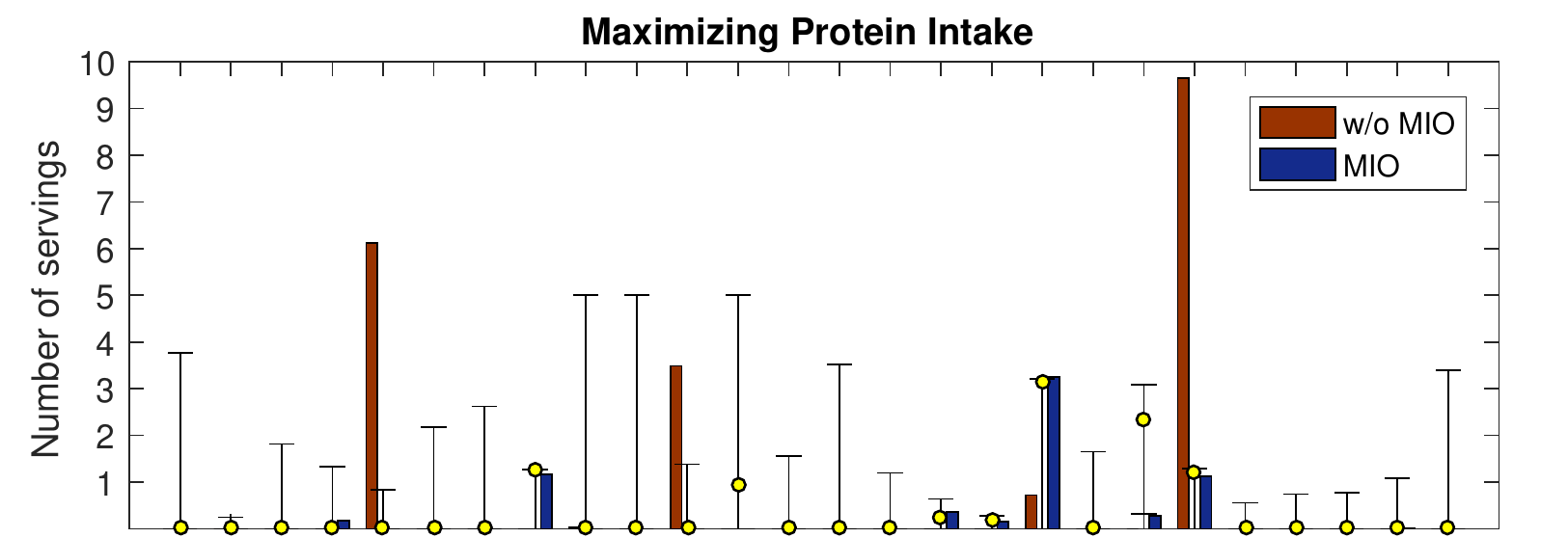}
    \includegraphics[width=.9\textwidth, trim= 0 5 0 0, clip=true]{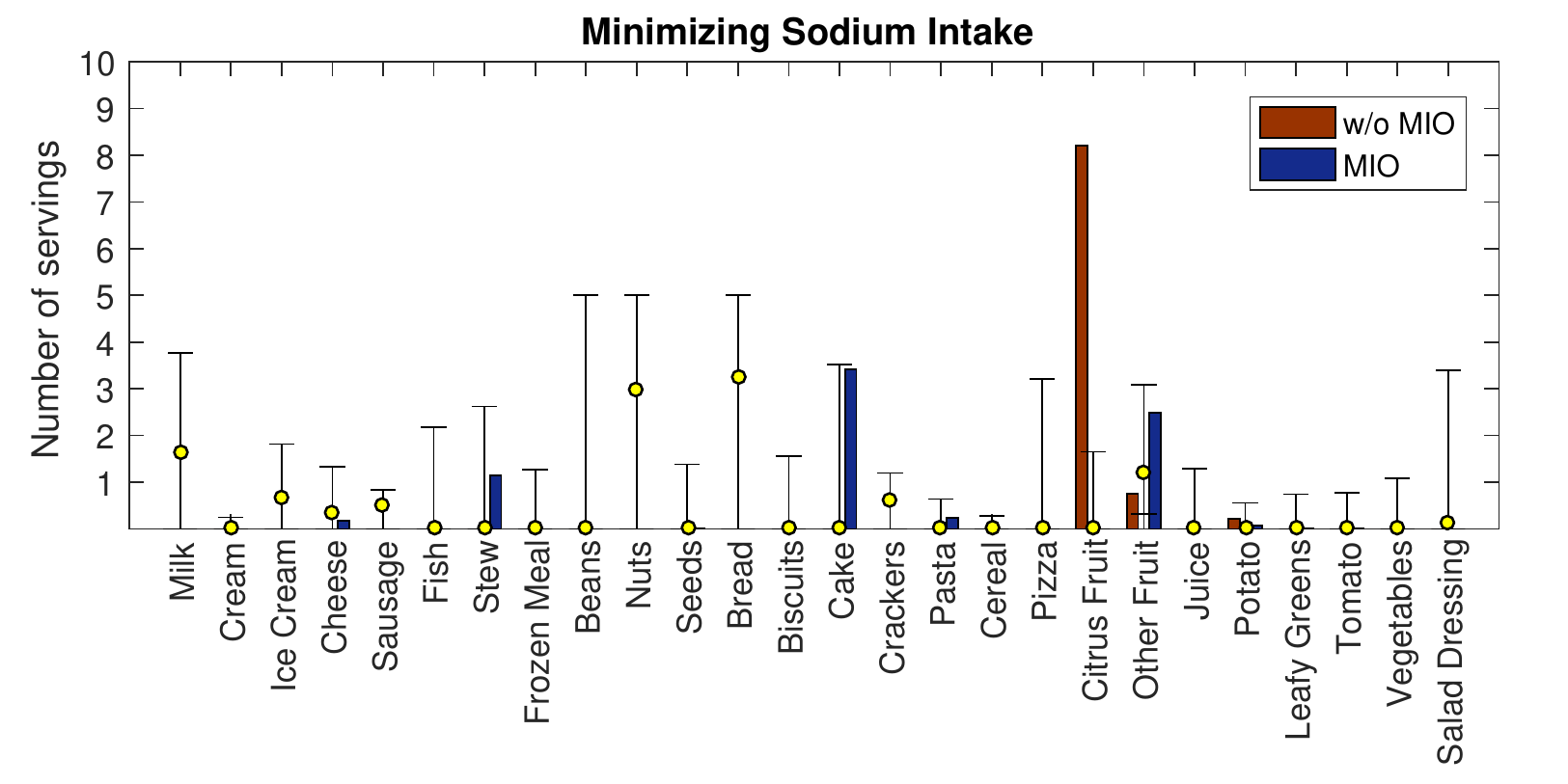}
    \caption{Comparison of diet recommendations with and without the imputed constraints. The range of past food consumption observations and the preferred observation are illustrated as error bars and yellow circles, respectively. 
    }\label{fig:barplot}
\end{figure}

\begin{figure}
\centering
\subfigure[\label{fig:spider-protein}]{   
    \includegraphics[width=0.47\textwidth, trim= 5 5 8 10, clip=true]{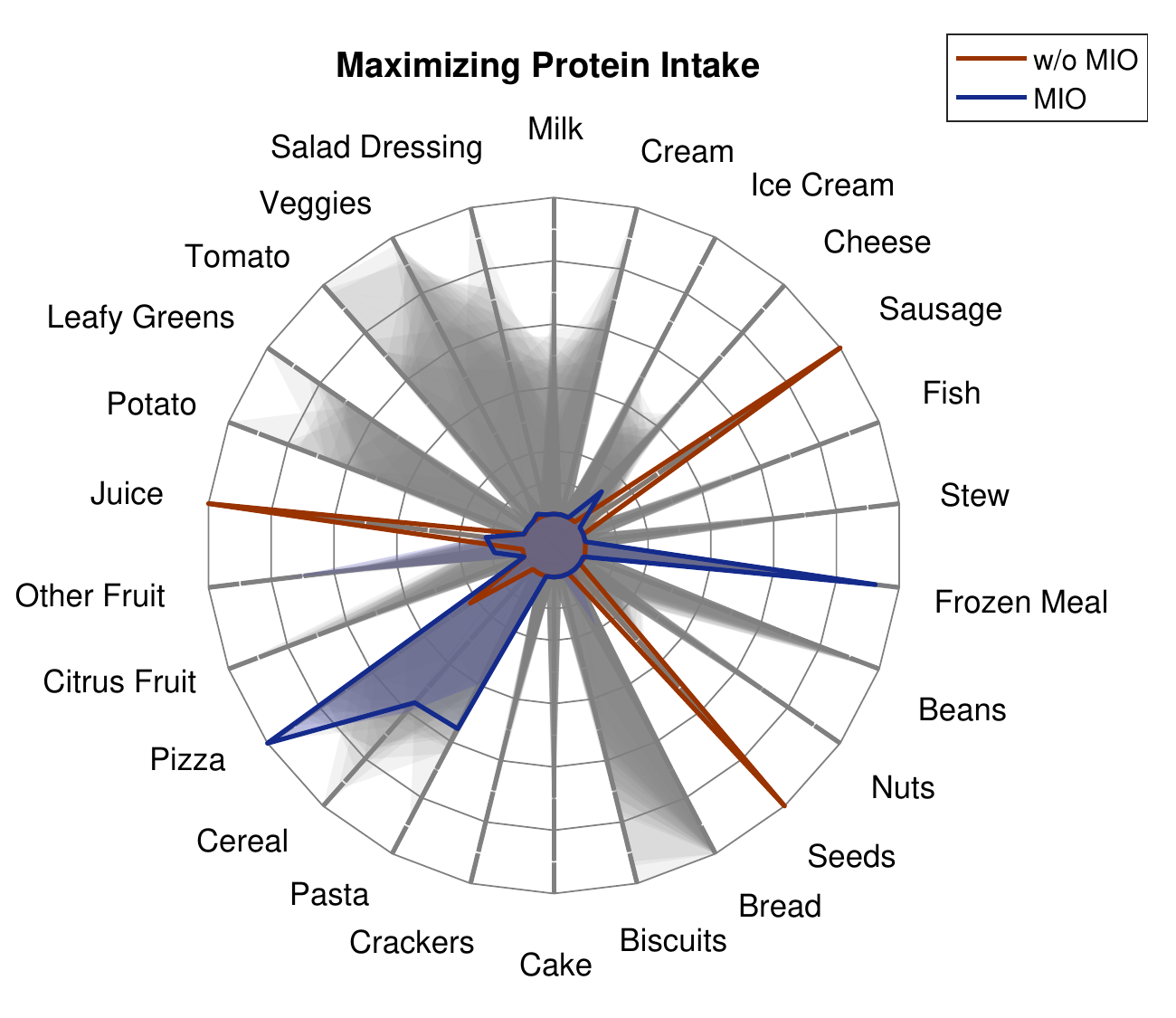}
}
\subfigure[\label{fig:spider-sodium}]{  
    \includegraphics[width=0.47\textwidth, trim= 5 5 8 10, clip=true]{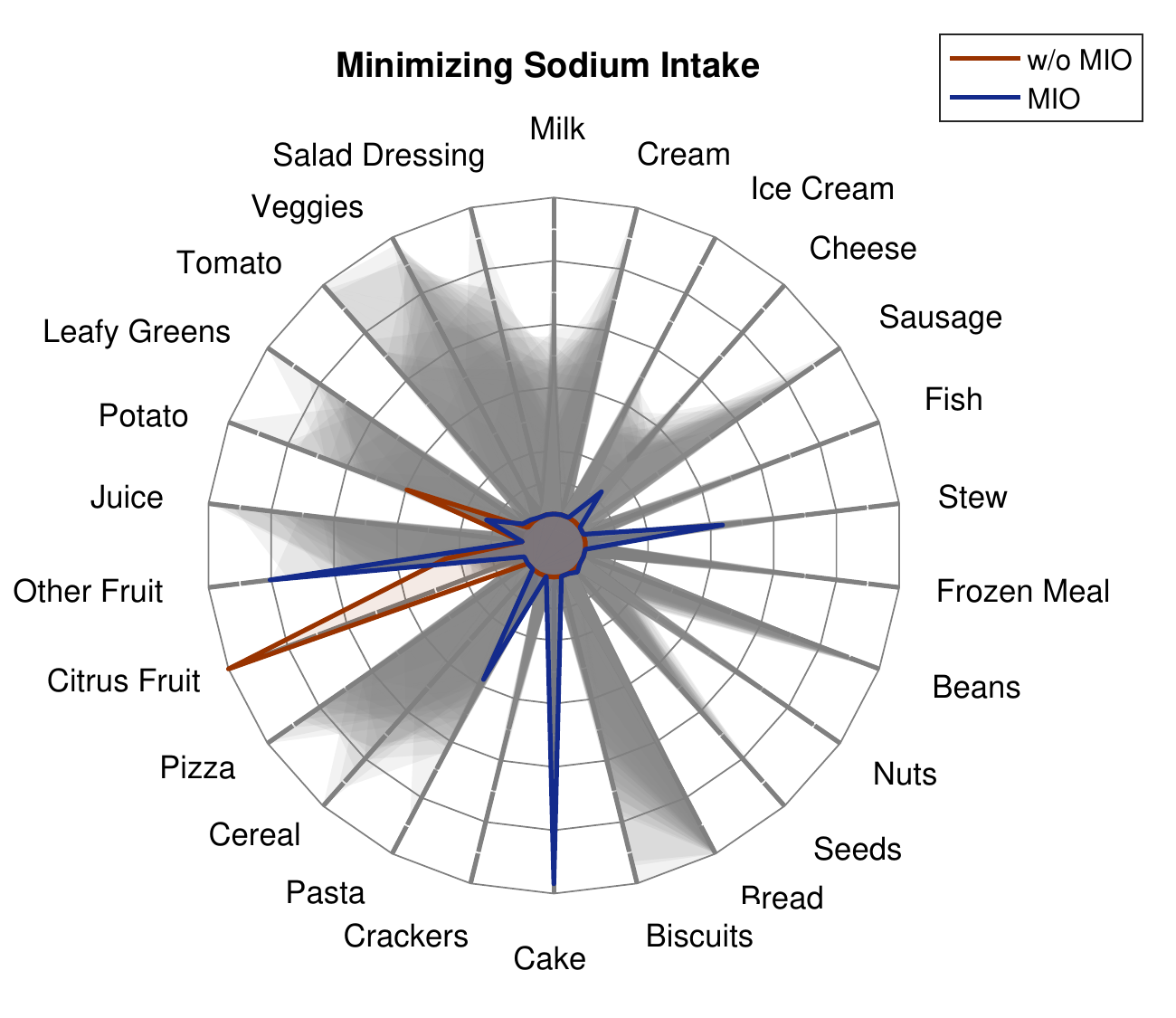}
}
\caption{Comparison of diet recommendations with and without the imputed constraints. The shaded grey areas show past observations. Darker grey shows a larger number of past food consumption for a food item.  }
\label{fig:spider}
\end{figure}

Figure~\ref{fig:barplot} shows the recommended diets for each of the two objectives with first considering only the known constraints (denoted as ``w/o MIO") and then also considering the additional imputed constraints (denoted as ``$\MIO$'') in red and blue bars, respectively. The ranges of the observations for each food item are shown by error bars, and the preferred observation ($\bx^0$) is highlighted in yellow circles on the error bar. Figure~\ref{fig:barplot} shows that the recommended diet without the imputed constraints includes larger quantities of fewer food items. For instance, when maximizing protein intake, over nine servings of juice and six servings of sausage are recommended to the user in the diet without $\MIO$ constraints. On the contrary, the diet with $\MIO$ constraints includes a variety of food items and more moderate quantities of each food item, closer to the food choices the user has made in the past. Similarly, when minimizing sodium intake,  without $\MIO$ constraints, the suggested diet only includes three food options, and the diet consists large amounts of fruit, while the diet with $\MIO$ constraints provides moderate amounts of six food items that more closely replicate meals consumed by the user in the past. 

To better visualize the different suggested diets in this multi-dimensional solution space, we also plotted these solutions using radar charts in Figure~\ref{fig:spider}. Each food item is individually Min-Max normalized. The diagrams serve to further illustrate the past data and compare each of the recommended diets on a relative scale, instead of the absolute scale showed in Figure~\ref{fig:barplot}. The gray shaded areas show the past observations where the darker areas show a larger amount of observed food intake. Similar to Figure~\ref{fig:barplot}, the solutions with and without $\MIO$ constraints are highlighted with blue and red lines, respectively. These plots confirm that in both cases, without $\MIO$ constraints, the diet often consists of food items that are not regularly consumed in the past and are limited in variety, while the diet with $\MIO$ constraints more closely resembles past observations. 

Finally, we quantify the differences between the diets recommended with and without the imputed $\MIO$ constraints for each objective by finding the average $L_1$~norm distance of all the observations from each of these diets as depicted in Table~\ref{tab:case3error}. In both cases, the diet with $\MIO$ constraints is closer to the observations. 

\begingroup
\setlength{\tabcolsep}{10pt} 
\renewcommand{\arraystretch}{1.2} 
\begin{table}[htbp]
    \centering
 \begin{tabular}{c c c} 
    \toprule
    & \multicolumn{2}{c}{Average $L_1$~norm distance} \\
    \cline{2-3}
Objective    & w/o $\MIO$ &  $\MIO$  \\
    \midrule
Maximizing Protein Intake   & 30.5 & 15.2 \\
Minimizing Sodium Intake    & 19.0 & 14.5 \\
    \bottomrule
\end{tabular}
    \caption{Average distance of recommended diets for each of the two objectives tested with and without the imputed constraints.}
    \label{tab:case3error}
\end{table}
\endgroup

\section{Model Extensions} \label{sec:discussions}
Typically, there are many assumptions in the inverse optimization models, depending on the particular structure of the model. In this section, we provide extensions for when some of these assumptions do not hold and discuss how our models can be adapted accordingly. 
In particular, we consider three extensions: (a) unknown cost vector, (b) noisy data, and (c) additional information on constraint parameters. We first briefly discuss the rationale behind the underlying standard assumptions and then provide extensions for lifting these assumptions.

\paragraph{(a) Unknown cost vector:}
When inferring the feasible region of a problem, the constraint parameters are unknown and hence, a set of given solutions cannot be labeled as feasible or infeasible. However, when a cost vector is known, the quality of the observation can be compared based on this cost vector. When, on the contrary, the cost vector is also unknown, there is no information about the quality of the given solutions. That is, we do not know whether a solution is feasible, and we are also not able to even assess or compare the quality of solutions. Such a problem setting may have limited practical use because it assumes experts have no knowledge about the objective or the constraints of the problem. Nevertheless, our models can be modified to consider an unknown cost vector. Let $\bc$ be a decision vector (unknown) and a candidate optimal solution $\bx^0$ be provided by experts. The updated model, denoted by $\MIO^{\bc}$, is
\begin{subequations}\label{eq:TMIOc}
\begin{align*}
\MIO^{\bc}: \underset{\bA, \bb, \by, \bw, \bc}{\text{minimize}} \quad & \mf(\bA, \bb; \bP), \\
\text{subject to} \quad & \eqref{eq:TMIO-PrimalFeas}-\eqref{eq:MIO-signs}, \\
\quad & \bc \in \mathbb{R}^n .
\end{align*}
\end{subequations} 
\noindent We note that the complexity of $\MIO^{\bc}$ is similar to that of $\MIO$ since the addition of $\bc$ as a variable does not introduce any new nonlinear terms into the model. The only difference is that $\bc$ is now a variable in the strong duality and the dual feasibility constraints~\eqref{eq:TMIO-StrongDual} and~\eqref{eq:TMIO-DualFeas}. 
%

\paragraph{(b) Noisy Data:} A standard assumption in many inverse optimization models in the literature is that the data is observed without noise. A few recent studies have considered that such perfect information may not be available, and even if the data is accurate, the models may be prone to overfitting to the given observations. In particular, the inverse model would always make $\bx^0$ exactly optimal for the forward problem and make any other solution that dominates $\bx^0$ infeasible.

To consider noisy data and avoid overfitting in our proposed model, a robust optimization approach such as that of \citet{ghobadi2018robust} can be employed by considering uncertainty sets around the observations. Such uncertainty sets can be considered around the preferred solution $\bx^0$ only, or alternatively around all observations $\bx^k, k \in \cK$. First, let $\cU^0$ be an uncertainty set around $\bx^0$. Since $\bc$ is known, the preferred solution within the uncertainty set $\cU^0$ is $\tilde{\bx}^0 =  \underset{\bx \in \cU^0}{\min} \{\bc' \bx^0\}$. Using this $\bx^0$ in the $\MIO$ formulation guarantees that the imputed feasible region is robust for all $\bx^0 \in \cU^0$. Next, assume we consider uncertainty sets around all observations. In this case, the feasibility of the uncertain observations needs to be guaranteed as well. A similar robust optimization approach can be employed to consider uncertainty sets $\cU^k$ around observations $\bx^k, \, \forall k \in \cK$. In addition to the strong duality constraint, in this case, the primal feasibility constraints are also modified to hold for any realization of the uncertainty set (\ie, $\forall \bx^k \in \cU^k$). 

\paragraph{(c) Additional Information on Constraint Parameters} 
In some inverse optimization models in the literature, some additional information (or side constraints) on the parameters of the forward model is considered in the inverse setting. Often, the purpose of this additional information is to avoid finding trivial (all-zero) solutions in the inverse model. Our inverse models avoid such trivial solutions by design and do not require the user to know and input such information on the inferred parameters. However, if such information exists, it can easily be incorporated into the model. Recall that we assumed that each constraint is either entirely known or entirely unknown. If additional information on some constraint parameters is available, the constraint (or set of constraints) is \emph{partially} known. In this case, we can modify the inverse model accordingly by assuming these partially-known constraints as part of the unknown constraints and adding the partial information as known constraints in the inverse model. For instance, if specific properties about the $i^{\text{th}}$ constraint is known (\eg, $b_i = \beta_i$ for $\beta_i \in \mathbb{R}$, or $a_{i1}\leq 2\,a_{i2}$), these properties can be explicitly added to the inverse model as known constraints. In general, if $\cA \subseteq \mathbb{R}^{m_1\times n}$ and $\cB \subseteq {R}^{m_1}$ capture the additional information on the LHS and RHS parameters, respectively, the $\MIO$ formulation can be adjusted by replacing constraint~\eqref{eq:MIO-signs} with $\bA \in \cA, \bb \in \cB.$
}

\section{Conclusions} \label{sec:conclusions}

Using inverse optimization to recover the feasible region can be applied to settings in which a set of solutions have previously been identified by experts as ``feasible" solutions, but the logic behind such labeling is not known. This paper proposes an inverse optimization approach for imputing fully- or partially-unknown constraint parameters of a forward optimization problem. The goal is to infer the feasible region of the forward problem such that all given observations become feasible and the preferred observation(s) become optimal. Identifying such feasible regions would provide a baseline for initial filtering of future observation as feasible or infeasible, before seeking an expert opinion. This information will, in turn, improve the flow of processes in expert-driven systems and reduce the time spent in manually identifying feasible solutions. In some applications, having such a data-driven approach would standardize the practice of quality control across different experts or different institutions.

We demonstrate the theoretical properties of our methodology and propose a new tractable reformulation for the nonlinear non-convex inverse model. We also present and discuss several loss functions that can be used to derive imputed feasible sets that have certain desired properties. Our numerical examples demonstrate the differences between these loss functions and serve as a basic guideline for users to choose the appropriate loss function  depending on the available data and the relevant application. \add{We further apply our methodology to a diet recommendation problem and show that the proposed model can impute the implicit constraints for each dieter and result in diet recommendations that are more palatable.}
An important future direction is to apply this methodology to a real-world large-scale dataset and to demonstrate the computational benefits of the proposed tractable reformulation methodology that allows for a more efficient solution of the originally nonlinear non-convex inverse optimization problems.

\clearpage
\bibliography{robustinverse2}

\begin{thebibliography}{38}
\expandafter\ifx\csname natexlab\endcsname\relax\def\natexlab#1{#1}\fi
\expandafter\ifx\csname url\endcsname\relax
  \def\url#1{\texttt{#1}}\fi
\expandafter\ifx\csname urlprefix\endcsname\relax\def\urlprefix{URL }\fi

\bibitem[{Ahuja and Orlin(2001)}]{Ahuja01}
Ahuja, R.~K., Orlin, J.~B., 2001. Inverse optimization. Operations Research
  49~(5), 771--783.

\bibitem[{Aswani et~al.(2018)Aswani, Shen, and Siddiq}]{Aswani16}
Aswani, A., Shen, Z.~J., Siddiq, A., 2018. Inverse optimization with noisy
  data. Operations Research 66~(3), 870--892.

\bibitem[{Ayer(2015)}]{Ayer14}
Ayer, T., 2015. Inverse optimization for assessing emerging technologies in
  breast cancer screening. Annals of Operations Research 230, 57--85.

\bibitem[{Babier et~al.(2018)Babier, Boutilier, Sharpe, McNiven, and
  Chan}]{babier2018inverse}
Babier, A., Boutilier, J.~J., Sharpe, M.~B., McNiven, A.~L., Chan, T.~C., 2018.
  Inverse optimization of objective function weights for treatment planning
  using clinical dose-volume histograms. Physics in Medicine \& Biology
  63~(10), 105004.

\bibitem[{Babier et~al.(2019)Babier, Chan, Lee, Mahmood, and
  Terekhov}]{babier2019ensemble}
Babier, A., Chan, T.~C., Lee, T., Mahmood, R., Terekhov, D., 2019. An ensemble
  learning framework for model fitting and evaluation in inverse linear
  optimization. arXiv preprint arXiv:1804.04576.

\bibitem[{Bertsimas et~al.(2012)Bertsimas, Gupta, and
  Paschalidis}]{Bertsimas12}
Bertsimas, D., Gupta, V., Paschalidis, I.~C., 2012. Inverse optimization: A new
  perspective on the {Black-Litterman} model. Operations Research 60~(6),
  1389--1403.

\bibitem[{Bertsimas et~al.(2015)Bertsimas, Gupta, and
  Paschalidis}]{Bertsimas15}
Bertsimas, D., Gupta, V., Paschalidis, I.~C., 2015. Data-driven estimation in
  equilibrium using inverse optimization. Mathematical Programming 153~(2),
  595--633.

\bibitem[{Birge et~al.(2017)Birge, Horta{\c{c}}su, and
  Pavlin}]{birge2017inverse}
Birge, J.~R., Horta{\c{c}}su, A., Pavlin, J.~M., 2017. Inverse optimization for
  the recovery of market structure from market outcomes: An application to the
  {MISO} electricity market. Operations Research 65~(4), 837--855.

\bibitem[{Brucker and Shakhlevich(2009)}]{brucker2009inverse}
Brucker, P., Shakhlevich, N.~V., 2009. Inverse scheduling with maximum lateness
  objective. Journal of Scheduling 12~(5), 475--488.

\bibitem[{{\v{C}}ern{\`y} and Hlad{\'\i}k(2016)}]{cerny2016inverse}
{\v{C}}ern{\`y}, M., Hlad{\'\i}k, M., 2016. Inverse optimization: towards the
  optimal parameter set of inverse {LP} with interval coefficients. Central
  European Journal of Operations Research 24~(3), 747--762.

\bibitem[{Chan et~al.(2019{\natexlab{a}})Chan, Eberg, Forster, Holloway,
  Ieraci, Shalaby, and Yousefi}]{chan2019inverse}
Chan, T.~C., Eberg, M., Forster, K., Holloway, C., Ieraci, L., Shalaby, Y.,
  Yousefi, N., 2019{\natexlab{a}}. An inverse optimization approach to
  measuring clinical pathway concordance. arXiv preprint arXiv:1906.02636.

\bibitem[{Chan and Kaw(2020)}]{chan2018inverse}
Chan, T.~C., Kaw, N., 2020. Inverse optimization for the recovery of constraint
  parameters. European Journal of Operational Research 282~(2), 415--427.

\bibitem[{Chan et~al.(2019{\natexlab{b}})Chan, Lee, and Terekhov}]{Chan15}
Chan, T.~C., Lee, T., Terekhov, D., 2019{\natexlab{b}}. Inverse optimization:
  Closed-form solutions, geometry, and goodness of fit. Management Science
  65~(3), 1115--1135.

\bibitem[{Chan et~al.(2014)Chan, Craig, Lee, and Sharpe}]{Chan14}
Chan, T. C.~Y., Craig, T., Lee, T., Sharpe, M.~B., 2014. Generalized inverse
  multi-objective optimization with application to cancer therapy. Operations
  Research 62~(3), 680--695.

\bibitem[{Chow and Recker(2012)}]{Chow12}
Chow, J. Y.~J., Recker, W.~W., 2012. Inverse optimization with endogenous
  arrival time constraints to calibrate the household activity pattern problem.
  Transportation Research Part B: Methodological 46~(3), 463--479.

\bibitem[{CSSE(2020)}]{CSSEDietData}
CSSE, 2020. {Diet Recommendation Data Page}.
  {http://systems.jhu.edu/data/health-care/diet} [Accessed: 19.06.2020].

\bibitem[{Dempe and Lohse(2006)}]{dempe2006inverse}
Dempe, S., Lohse, S., 2006. Inverse linear programming. In: Recent Advances in
  Optimization. Springer, pp. 19--28.

\bibitem[{Dong and Zeng(2018)}]{dong2018inferring}
Dong, C., Zeng, B., 2018. Inferring parameters through inverse multiobjective
  optimization. arXiv preprint arXiv:1808.00935.

\bibitem[{Erkin et~al.(2010)Erkin, Bailey, Maillart, Schaefer, and
  Roberts}]{Erkin10}
Erkin, Z., Bailey, M.~D., Maillart, L.~M., Schaefer, A.~J., Roberts, M.~S.,
  2010. Eliciting patients' revealed preferences: An inverse {M}arkov decision
  process approach. Decision Analysis 7~(4), 358--365.

\bibitem[{Esfahani et~al.(2018)Esfahani, Shafieezadeh-Abadeh, Hanasusanto, and
  Kuhn}]{esfahani2018IncompleteInfo}
Esfahani, P.~M., Shafieezadeh-Abadeh, S., Hanasusanto, G.~A., Kuhn, D., 2018.
  Data-driven inverse optimization with incomplete information. Mathematical
  Programming 167~(1), 191--234.

\bibitem[{Fourer et~al.(1993)Fourer, Gay, and Kernighan}]{ampl}
Fourer, R., Gay, D., Kernighan, B., 1993. {AMPL}. Boyd and Fraser.

\bibitem[{Ghate(2020)}]{ghate2020inverse}
Ghate, A., 2020. Inverse optimization in semi-infinite linear programs.
  Operations Research Letters 48~(3), 278--285.

\bibitem[{Ghobadi et~al.(2018)Ghobadi, Lee, Mahmoudzadeh, and
  Terekhov}]{ghobadi2018robust}
Ghobadi, K., Lee, T., Mahmoudzadeh, H., Terekhov, D., 2018. Robust inverse
  optimization. Operations Research Letters 46~(3), 339--344.

\bibitem[{Goli et~al.(2018)Goli, Boutilier, Craig, Sharpe, and Chan}]{Goli15}
Goli, A., Boutilier, J.~J., Craig, T., Sharpe, M.~B., Chan, T.~C., 2018. A
  small number of objective function weight vectors is sufficient for automated
  treatment planning in prostate cancer. Physics in Medicine \& Biology
  63~(19), 195004.

\bibitem[{G{\"u}ler and Hamacher(2010)}]{guler2010capacity}
G{\"u}ler, {\c{C}}., Hamacher, H.~W., 2010. Capacity inverse minimum cost flow
  problem. Journal of Combinatorial Optimization 19~(1), 43--59.

\bibitem[{Iyengar and Kang(2005)}]{Iyengar05}
Iyengar, G., Kang, W., 2005. Inverse conic programming with applications.
  Operations Research Letters 33~(3), 319--330.

\bibitem[{Keshavarz et~al.(2011)Keshavarz, Wang, and Boyd}]{Keshavarz11}
Keshavarz, A., Wang, Y., Boyd, S., 2011. Imputing a convex objective function.
  In: 2011 {IEEE} International Symposium on Intelligent Control (ISIC). IEEE,
  pp. 613--619.

\bibitem[{Naghavi et~al.(2019)Naghavi, Foroughi, and
  Zarepisheh}]{naghavi2019inverse}
Naghavi, M., Foroughi, A.~A., Zarepisheh, M., 2019. Inverse optimization for
  multi-objective linear programming. Optimization Letters 13~(2), 281--294.

\bibitem[{Petit and Trapp(2019)}]{cplex}
Petit, T., Trapp, A.~C., 2019. Enriching solutions to combinatorial problems
  via solution engineering. INFORMS Journal on Computing 31~(3), 429--444.

\bibitem[{Saez-Gallego and Morales(2018)}]{saez2018short}
Saez-Gallego, J., Morales, J.~M., 2018. Short-term forecasting of
  price-responsive loads using inverse optimization. IEEE Transactions on Smart
  Grid 9~(5), 4805--4814.

\bibitem[{Saunders and Murtagh(2003)}]{saunders2003minos}
Saunders, M.~A., Murtagh, B.~A., 2003. {MINOS} 5.51 user's guide.

\bibitem[{Shahmoradi and Lee(2019)}]{shahmoradi2019quantile}
Shahmoradi, Z., Lee, T., 2019. Quantile inverse optimization: Improving
  stability in inverse linear programming. arXiv preprint arXiv:1908.02376.

\bibitem[{Tavasl{\i}o{\u{g}}lu et~al.(2018)Tavasl{\i}o{\u{g}}lu, Lee, Valeva,
  and Schaefer}]{tavasliouglu2018structure}
Tavasl{\i}o{\u{g}}lu, O., Lee, T., Valeva, S., Schaefer, A.~J., 2018. On the
  structure of the inverse-feasible region of a linear program. Operations
  Research Letters 46~(1), 147--152.

\bibitem[{Troutt et~al.(2008)Troutt, Brandyberry, Sohn, and
  Tadisina}]{Troutt08}
Troutt, M.~D., Brandyberry, A.~A., Sohn, C., Tadisina, S.~K., 2008. Linear
  programming system identification: The general nonnegative parameters case.
  European Journal of Operational Research 185~(1), 63--75.

\bibitem[{Troutt et~al.(2006)Troutt, Pang, and Hou}]{Troutt06}
Troutt, M.~D., Pang, W.~K., Hou, S.~H., 2006. Behavioral estimation of
  mathematical programming objective function coefficients. Management Science
  52~(3), 422--434.

\bibitem[{{U.S. Department of Health and Human Services and U.S. Department of
  Agriculture}(2015)}]{DietData_Const}
{U.S. Department of Health and Human Services and U.S. Department of
  Agriculture}, 2015. {2015--2020 Dietary Guidelines for Americans, 8th
  Edition}.
  {https://health.gov/our-work/food-nutrition/2015-2020-dietary-guidelines/guidelines/appendix-7/\#table-a7-1}
  [Accessed: 19.06.2020].

\bibitem[{Zhang and Liu(1996)}]{zhang1996calculating}
Zhang, J., Liu, Z., 1996. Calculating some inverse linear programming problems.
  Journal of Computational and Applied Mathematics 72~(2), 261--273.

\bibitem[{Zhang and Liu(1999)}]{zhang1999further}
Zhang, J., Liu, Z., 1999. A further study on inverse linear programming
  problems. Journal of Computational and Applied Mathematics 106~(2), 345--359.

\end{thebibliography}
\bibliographystyle{elsarticle-harv}

\clearpage
{\color{myGreen}

\APPENDIX{Data Summary} 
\label{appendix}
Table~\ref{tab:case3_fooddetails} shows the summary of the 100 observations (\ie, days) of daily food intake.  The column `count' shows how many times each food was consumed (\ie, in how many of the 100 observations). The next two columns show the average and standard deviation of the number of servings of each food, in the days that the food was consumed. 

\begin{table}[htbp]
\small
    \centering
 \begin{tabular}{l c c c}
    \toprule
    Food Item 	&	Count	&	Avg. Consumption	&	Std. Dev.	\\\midrule
Milk	&	60	&	1.4	&	1.0	\\
Cream	&	20	&	0.2	&	0.0	\\
Ice Cream	&	40	&	0.9	&	0.3	\\
Cheese	&	40	&	0.6	&	0.3	\\
Sausage	&	20	&	0.6	&	0.1	\\
Fish	&	20	&	1.5	&	0.4	\\
Stew	&	20	&	2.1	&	0.4	\\
Frozen Meal	&	40	&	0.8	&	0.2	\\
Beans	&	20	&	5.0	&	0.1	\\
Nuts	&	60	&	1.8	&	1.7	\\
Seeds	&	20	&	0.9	&	0.2	\\
Bread	&	80	&	2.9	&	1.5	\\
Biscuits	&	20	&	1.0	&	0.3	\\
Cake	&	20	&	2.7	&	0.5	\\
Crackers	&	40	&	0.6	&	0.2	\\
Pasta	&	40	&	0.4	&	0.1	\\
Cereal	&	20	&	0.2	&	0.0	\\
Pizza	&	20	&	2.5	&	0.4	\\
Citrus Fruit	&	20	&	1.1	&	0.3	\\
Other Fruit	&	100	&	1.5	&	0.6	\\
Juice	&	20	&	0.9	&	0.2	\\
Potato	&	20	&	0.3	&	0.1	\\
Leafy Greens	&	20	&	0.5	&	0.1	\\
Tomato	&	20	&	0.5	&	0.2	\\
Vegetables	&	40	&	0.5	&	0.4	\\
Salad Dressing	&	40	&	1.1	&	1.0	\\
    \bottomrule
\end{tabular}
    \caption{For every consumed food item, the number of consumption, the average serving size, and the standard deviation is provided in columns Count, Avg. Consumption, and Std. Dev., respectively. }
    \label{tab:case3_fooddetails}
\end{table}
}

\end{document}